\documentclass[reqno]{amsart}[standalone]
\usepackage[T1]{fontenc}
\usepackage{amsmath}
\usepackage{amssymb}
\usepackage{amsthm}
\usepackage{amsfonts}
\usepackage{url}
\usepackage{comment}
\usepackage{graphicx}
\usepackage{bm}
\usepackage{glossaries}
\usepackage{algpseudocode}
\usepackage[ruled,linesnumbered]{algorithm2e}
\usepackage{mathtools}
\usepackage[inline]{enumitem}
\usepackage{cases}
\usepackage{hyperref}
\usepackage{tikz}
\usetikzlibrary{arrows.meta, positioning, calc}
\usepackage{subcaption}
\usepackage{xcolor}
\usepackage{float}
\usepackage{siunitx}
\usepackage[top=1in, bottom=1.25in, left=1.25in, right=1.25in]{geometry}

\usepackage{mathrsfs}
\usepackage[foot]{amsaddr}
\usepackage{caption}
\usepackage{mathbbol}
\usepackage{biblatex}
\graphicspath{{Images/}}

\theoremstyle{definition}
\newtheorem{definition}{Definition}[section]
\theoremstyle{remark}
\newtheorem{remark}[definition]{Remark}

\newcommand{\R}{\mathbb{R}}

\newcommand{\bs}{\boldsymbol}
\DeclareMathOperator*{\argmin}{\arg\,\min}

\newacronym{nn}{NN}{Neural Network}
\newacronym{mor}{MOR}{Model Order Reduction}
\newacronym{rom}{ROM}{Reduced Order Model}
\newacronym{pod}{POD}{Proper Orthogonal Decomposition}
\newacronym{sindy}{SINDy}{Sparse Identification of Nonlinear Dynamics}
\newacronym{cfd}{CFD}{Computational Fluid Dynamics}
\newacronym{ae}{AE}{Autoencoder}
\newacronym{fe}{FE}{Finite Element}
\newacronym{fv}{FV}{Finite Volume}
\newacronym{fom}{FOM}{full-order method}
\newacronym{stlsq}{STLSQ}{Sequential Thresholded Least-Squares}

\begin{document}
\title[]{Sparse Identification for bifurcating phenomena in Computational Fluid Dynamics}
\author{Lorenzo Tomada$^{1}$, Moaad Khamlich$^{1}$, Federico Pichi$^{1}$, Gianluigi Rozza$^{1}$}
\address{$^1$ mathLab, Mathematics Area, SISSA, via Bonomea 265, I-34136 Trieste, Italy}
\email{\texttt{\{ltomada, mkhamlic, fpichi, grozza\}@sissa.it}}
\newif\ifincludeKeywords
\includeKeywordstrue 

\begin{abstract}
This work investigates model reduction techniques for nonlinear parameterized and time-dependent PDEs, specifically focusing on bifurcating phenomena in \gls{cfd}.
We develop interpretable and non-intrusive \glspl{rom} capable of capturing dynamics associated with bifurcations by identifying a minimal set of coordinates.

Our methodology combines the \gls{sindy} method with a deep learning framework based on \gls{ae} architectures.
To enhance dimensionality reduction, we integrate a nested \gls{pod} with the SINDy-AE architecture, enabling a sparse discovery of system dynamics while maintaining efficiency of the reduced model.

We demonstrate our approach via two challenging test cases defined on sudden-expansion channel geometries: a symmetry-breaking bifurcation and a Hopf bifurcation.
Starting from a comprehensive analysis of their high-fidelity behavior, i.e.\ symmetry-breaking phenomena and the rise of unsteady periodic solutions, we validate the accuracy and computational efficiency of our \glspl{rom}. 

The results show successful reconstruction of the bifurcations, accurate prediction of system evolution for unseen parameter values, and significant speed-up compared to full-order methods.

\ifincludeKeywords

\fi
\end{abstract}

\maketitle
\glsresetall

\section{Introduction and Motivation}
\label{sec:intro}

An important topic in \gls{cfd} is the study of flow transitions in fluids and gases resulting from changes in parameters.
Several phenomena are related to the loss of stability and uniqueness for flow configurations, such as droplet formation in capillary flows~\cite{Dijkstra_Wubs_2023, Dijkstra_Wubs_Cliffe_Doedel_Dragomirescu_Eckhardt_Gelfgat_Hazel_Lucarini_Salinger_et_al_2014}, or transition from laminar to unsteady regime in flows past a cylinder \cite{Zdravkovich1997}.
The study of such qualitative changes extends to various applications, including pipe systems \cite{Gavilan22}, branching networks in biological systems and engineering \cite{PRADHAN2019108483}, and clinical scenarios such as cardiovascular system modeling \cite{Pitton_2017_2}.
While these physical processes may vary significantly, they often share underlying mathematical structures. 

Bifurcation theory provides a valuable framework for understanding fluid flow instabilities at an abstract level, offering a structured approach to evaluate how solutions of mathematical models respond to variations in the parameter configuration \cite{Dijkstra_Wubs_2023}.

The objective of this work is to investigate complex bifurcating phenomena in \gls{cfd} using efficient \glspl{rom} \cite{RozzaBallarinScandurraPichi2024, Benner2017, MOR}.
We consider two relevant test cases, each representing a different type of bifurcation, involving both steady-state and periodic unsteady equilibrium configurations.

The first test case examines the \emph{Coandă effect} in a two-dimensional sudden-expansion channel, where the flow exhibits asymmetric behavior and wall attachment \cite{PichiDrivingBifurcatingParametrized2022a,Khamlich_2022}.
This system undergoes a \emph{supercritical pitchfork bifurcation}, with solution multiplicity and stability depending on the Reynolds number $Re$.
This setting serves as a simplified model for a \textit{mitral valve} geometry to study a heart condition characterized by abnormal blood flow between the left ventricle and the left atrium due to impaired mitral valve closure \cite{tritton1977physical}. Indeed, this framework effectively captures essential properties of blood flow between heart chambers, where wall-hugging behavior can lead to inaccurate echocardiographic measurements.
The symmetry breaking in the sudden-expansion channel, despite geometric and boundary conditions symmetry, is due to the presence of the nonlinear term, causing the ill-posedness of the parametric problem and the non-uniqueness of the solution's behavior. 

The second test case investigates a benchmark problem in \gls{cfd}, focusing on identifying and reconstructing the \textit{Hopf bifurcation} that marks the transition from steady to unsteady behavior, where periodic solutions emerge \cite{QUAINI}.
This bifurcation occurs at higher Reynolds numbers than the symmetry-breaking one and is driven by the loss of stability of the steady solution, originating a family of unsteady solutions that branches from it.
We chose the sudden-expansion channel framework due to its rich dynamics, exhibiting multiple bifurcations even at low Reynolds numbers due to the Navier-Stokes equations' nonlinearity.

The full-order numerical approaches usually exploited to reconstruct the flow dynamics employ the \gls{fe} method \cite{Quarteroni_2017} and the \gls{fv} method \cite{FVM}.

Despite the extensive research on bifurcations in \gls{cfd} and advances in scientific computing, accurate numerical simulations remain computationally intensive \cite{TESIFP}.
This challenge stems from solution non-uniqueness and the many-query context to understand critical points and post-bifurcation dynamics, and the computational burden increases even further in the considered time-dependent setting. 

To address these challenges and enable real-time analysis, we rely on reduced order modeling techniques \cite{hesthaven_certified_2015,quarteroni_reduced_2016,RozzaBallarinScandurraPichi2024}. In particular, this work explores \glspl{rom} based on \gls{sindy} \cite{Brunton_2016}, a non-intrusive method that explicitly derives equations describing the dynamical system from high-fidelity snapshots. To handle high-dimensional data, we incorporate deep learning techniques using \glspl{nn} \cite{Champion_2019, DEEPDELAY, Conti_2023}, combining \gls{sindy}'s parsimony and interpretability with the universal approximation capabilities of deep \glspl{nn}.

The main contributions of this work include: (i) the extension of previous research that only focused on steady-state problems \cite{TESIFP, Khamlich_2022, PICHI2023105813, gonnella2025stochasticperturbationapproachnonlinear}, by providing a comprehensive investigation of time-dependent bifurcating behavior for both benchmarks, and (ii) their efficient identification via a novel integration of \emph{nested} \gls{pod} \cite{nested_POD1, nested_POD2} with an \gls{ae}-based \gls{sindy} approach. The developed methodology is efficient, purely data-driven, and interpretable, with strong capabilities for temporal and parametric extrapolation.

Our results show that it is fundamental to perform a careful analysis of such non-trivial behaviors, already at full-order level, investigating potentially impactful symmetry-breaking phenomena that affect the transition from steady to unsteady equilibria, and the time needed to reach convergence.
Moreover, the development and analysis of such problems via \glspl{rom} allows an efficient identification of the latent evolution of the system, and the consequent reconstruction of full-order trajectories and bifurcation diagrams.

The work is organized as follows:
\begin{itemize}
    \item Section \ref{sec:related} reviews the state of the art regarding bifurcations in \gls{cfd}, \gls{rom} and \gls{sindy} approaches;
    \item Section \ref{sec_bifurcations} details the mathematical framework, and presents the bifurcating behaviors obtained via high-fidelity numerical methods;
    \item Section \ref{sec:ROM} introduces the reduced order modeling tools, i.e.\ \gls{sindy} and nested \gls{pod}, and their integration in the SINDy-AE-nested-POD architecture;
    \item Section \ref{sec:numerical} illustrates the numerical results for symmetry-breaking and Hopf bifurcations.
\end{itemize}

\section{Related Works}\label{sec:related}

This work investigates two distinct bifurcating phenomena in fluid flows in sudden-expansion channels: a symmetry-breaking bifurcation known as the Coandă effect \cite{PichiDrivingBifurcatingParametrized2022a,Khamlich_2022,Pitton_2017}, and a Hopf bifurcation \cite{QUAINI,Fortin_localization}. 
The Coandă effect arises from the tendency of the fluid jet to attach to nearby surfaces, driven by velocity fluctuations that generate transverse pressure gradients, which ultimately maintain flow asymmetry~\cite{tritton1977physical, wille1965coanda}. The physical mechanism has been extensively studied \cite{Cherdron_Durst_Whitelaw_1978, Sobey_Drazin_1986}, revealing that at low Reynolds numbers, viscous dissipation stabilizes the flow \cite{HAWA_RUSAK_2001}. As Reynolds number increases, convective effects upstream of the expansion overcome this stabilization, leading to symmetric solution instability.

The seminal work by Sobey and Drazin \cite{Sobey_Drazin_1986} combined asymptotic, numerical, and experimental approaches to identify a Hopf bifurcation, with the bifurcation point later refined by Quaini et al.\ \cite{QUAINI}. The emergence of periodic solutions through Hopf bifurcations for different benchmarks is well-documented in \gls{cfd} literature \cite{Sobey_Drazin_1986, Fortin_localization, ArioliKoch2021, DusekNumericalTheoreticalStudy1994}.

Recent years have seen significant developments in \glspl{rom} for these systems, extending both the range of applications and the available methodologies. To reduce the gap with real-world problems, more advanced physical models have been considered, e.g.\ including fluid-structure interaction effects \cite{Khamlich_2022}, optimization and control to drive bifurcating phenomena \cite{PichiDrivingBifurcatingParametrized2022a,Boull2023}, interpretable \glspl{rom} for the fluidic pinball \cite{Pastur_2019}, geometrically parametrized domains \cite{PICHI2023105813,Bravo_2024}, and the compressible \gls{cfd} setting \cite{Tonicello_2024}. As concerns the methodological advancements,  deflated approaches have been used to discover the multiple coexisting states \cite{pintore_2021}, recently combined with greedy algorithms for efficiency, adaptivity and model certification \cite{pichi2025deflationbasedcertifiedgreedyalgorithm}. A stochastic perturbation approach \cite{gonnella2025stochasticperturbationapproachnonlinear} has been developed to obtain insights about the bifurcation points in the uncertainty quantification setting. Moreover, localized \gls{rom} strategies have been used to build specific reduced spaces depending on the fluid's regime \cite{Hess_2019}.

From a different perspective, an interesting and potentially effective approach consists in tackling the problem by means of sparse identification, in particular employing \gls{sindy} \cite{Brunton_2016}. Originally developed for learning dynamics from system snapshots, \gls{sindy} has evolved through various generalizations \cite{SINDy_ensemble, SINDyCP, Messenger_2021} and found applications in fluid flow identification \cite{Fukami_Murata_Zhang_Fukagata_2021, wang_zhang_2024} and convective phenomena. The integration of \gls{sindy} with \glspl{ae} for dimensionality reduction was pioneered by Champion et al.~\cite{Champion_2019}, extended in~\cite{DEEPDELAY}, and adapted for parameterized systems by~\cite{Conti_2023}.

At the same time, neural networks and data-driven approaches have become increasingly prevalent in dynamical systems research, particularly for bifurcation analysis \cite{Kalia2021LearningNF,PICHI2023105813} and low-dimensional flow representation \cite{Della_Pia_2024, spider}. However, challenges persist regarding interpretability and overfitting tendencies. The \gls{nn}-\gls{sindy} approach was developed specifically to address these limitations.

Our work advances the state of the art by combining the SINDy-AE approach with nested \gls{pod} \cite{nested_POD1, nested_POD2}, a dimensionality reduction technique well-suited for time-dependent parameterized PDEs. This integration improves computational efficiency during the offline phase, enabling the processing of larger snapshot sets with reduced computational resources. While \gls{sindy} has been applied to PDE identification \cite{ZHENG2026114443, Rudy_2017, Messenger_2021}, we maintain an ODE framework, consistent with recent developments \cite{Champion_2019, Conti_2023, VVV}.

\section{Bifurcations in the Navier-Stokes equations}\label{sec_bifurcations}

We consider the strong form of the incompressible unsteady Navier-Stokes equations describing the behavior of fluid flows with constant density, given by\footnote{The arguments of the functions $\bs{u}, \bs{f}$, $p$ and their derivatives are omitted for notational simplicity.}:
\begin{numcases}{}
    \frac{\partial\bs{u}}{\partial t} + (\bs{u}\cdot\nabla)\bs{u} - \mu\Delta\bs{u} + \nabla p = \bs{f}, & $\forall(\bs{x}, t)\in\Omega\times\R^+$, \label{eq:NS1}\\
    \nabla\cdot \bs{u}= 0, & $\forall(\bs{x}, t)\in\Omega\times\R^+$, \label{eq:NS2}
\end{numcases}
where $\mu\in\R^+$ represents the kinematic viscosity, $\bs{f}\in L^2(\R^+; [L^2(\Omega)]^d)$ denotes the force per unit volume, and $\bs{u}$ and $p$ represent the velocity field and pressure of the fluid, respectively. The domain $\Omega\subseteq\mathbb{R}^d$ with $d = 2, 3$, and the nonlinear term $(\bs{u}\cdot\nabla)\bs{u}$ introduces convective transport.

In the next sections, we present and study the bifurcating phenomena arising from these equations when posed on sudden-expansion channel geometries: a symmetry-breaking bifurcation and a Hopf bifurcation.

\subsection{Symmetry-breaking bifurcation}\label{subsec:coanda_problem}
\begin{figure}[ht]
    \centering
    \includegraphics[width=0.9\textwidth]{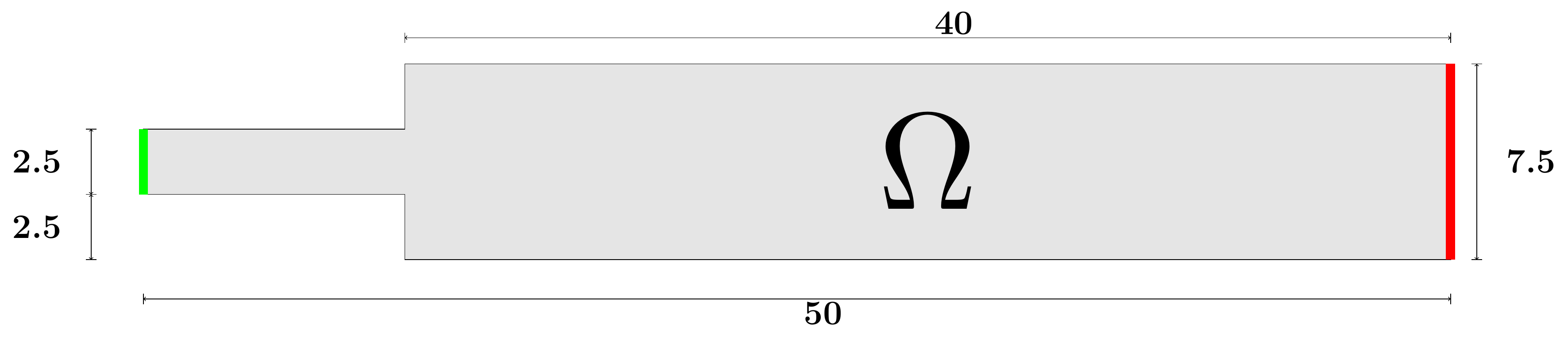}
    \caption{Problem geometry showing inlet boundary $\Gamma_\text{in}$ (green), outlet boundary $\Gamma_\text{out}$ (red), and wall boundary $\Gamma_0$ (black).}
    \label{figMesh}
\end{figure}

We consider the sudden-expansion channel $\Omega\subset\R^2$, depicted in Figure \ref{figMesh}, and defined as $\Omega = \bigl((0, 10)\times (2.5, 5)\bigr)\cup \bigl((10, 50)\times(0, 7.5)\bigr).$

We solve Equations \eqref{eq:NS1}-\eqref{eq:NS2} in $\Omega$ with parametric range $\mu\in\mathbb{P}=[0.5, 2]$, and subject to the following boundary conditions:
\begin{equation}\label{eq:bc_coanda}
    \begin{dcases}
    \bs{u}(\bs{x}, 0) = \bs{0}, &\text{in } \Omega, \\
    \bs{u}(\bs{x}, t) = [20(y - 2.5)(5 - y),\, 0], &\text{on } \Gamma_{\text{in}}\times(0,T), \\
    \bs{u}(\bs{x}, t) = \bs{0}, &\text{on } \Gamma_0\times(0,T), \\
    (-\mu\nabla\bs{u} + p\mathbf{I})\bs{n} = \bs{0}, &\text{on } \Gamma_{\text{out}}\times(0,T),
    \end{dcases}
\end{equation}
where $\Gamma_{\text{in}}=\{0\}\times [2.5,\,5]$, $\Gamma_{\text{out}} = \{50\}\times[0,\, 7.5]$, and $\Gamma_0 = \partial\Omega\smallsetminus(\Gamma_{\text{in}}\cup\Gamma_{\text{out}})$.

The Reynolds number, balancing convection and diffusion effects, is defined as $Re = UL/\mu$, where $L$ and $U$ are characteristic length and velocity for the fluid problem, respectively chosen as the inlet length $L=2.5$, and the maximum inlet velocity $U=31.25$. This way, considering the parametric setting $\mu\in[0.5,\,2]$, the Reynolds number range is given by $Re\in[39,\, 156]$.

For the numerical simulations, we employ the \gls{fe} method implemented in FEniCS \cite{LoggEtal_2012}. The spatial discretization uses a mesh with $2752$ triangular elements and Taylor--Hood $\mathbb{P}_2-\mathbb{P}_1$ pair for velocity-pressure fields, ensuring the validity of the \textit{inf-sup} condition \cite{Rozza_2007} while maintaining computational efficiency. 

For time integration, we use the $\theta$-method with $\theta = \frac{1}{2}$ (trapezoidal rule), which provides second-order accuracy and unconditional stability, with a time step $dt = 10^{-2}$. The nonlinear system at each time step is solved using Newton iterations with a direct LU solver for the linear steps. 
Finally, we employ the following relative convergence criterion between successive iterations:
\begin{equation}\label{eqSTOP}
\frac{\|\bs{u}^{n+1}-\bs{u}^n\|}{\|\bs{u}^n\|}<\tau,
\end{equation}
where $\tau=10^{-7}$ is chosen to ensure proper resolution of the bifurcation behavior.

\subsubsection{Bifurcation Analysis}
Previous studies have revealed a rich bifurcation structure for the steady case \cite{QUAINI, pintore_2021}. Indeed, as described in Section \ref{sec:intro}, there exists a critical viscosity $\mu^\ast\approx 0.96$ such that the system undergoes a supercritical pitchfork bifurcation \cite{PRACTICALBIF}, and the uniqueness of the solution is lost.
For lower viscosity values, the symmetric solution becomes unstable, the inertial effects begin to dominate, and two stable asymmetric solutions emerge, with the jet bending toward either the upper or lower wall. 
To comprehensively study the transition from uniqueness to bifurcating regions, we choose the parametric space $\mathbb{P}=[0.5, 2]$ so that it includes the bifurcation point $\mu^\ast$.

\begin{figure}[ht]
\centering
\begin{subfigure}{0.49\linewidth}
\centering
\includegraphics[width=\linewidth, trim=0 10 0 130, clip]{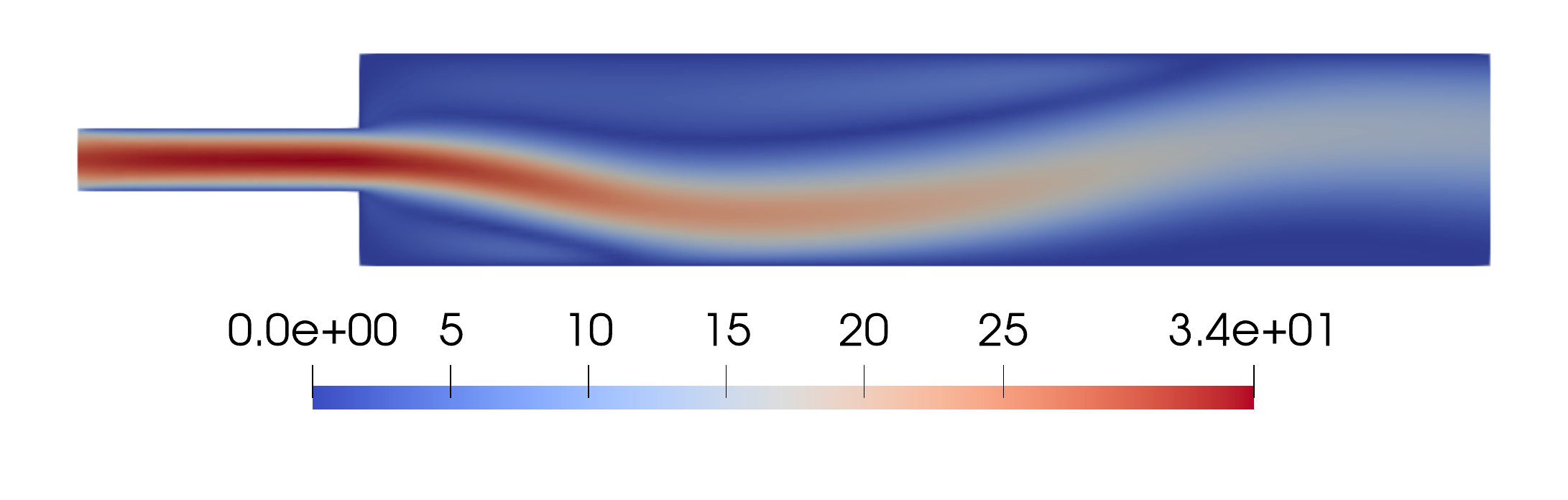}
\caption{Stable asymmetric solution}
\label{fig:stable_0.5}
\end{subfigure}\hfill
\begin{subfigure}{0.49\linewidth}
\centering
\includegraphics[width=\linewidth, trim=0 10 0 130, clip]{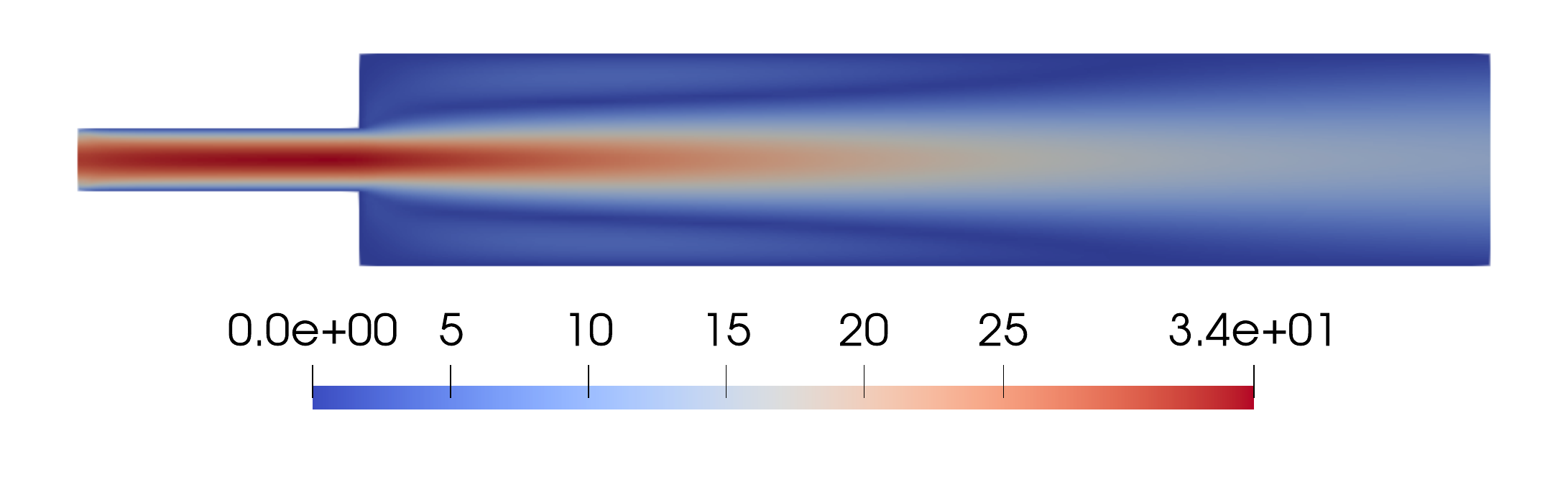}
\caption{Unstable symmetric solution}
\label{fig:unstable_0.5}
\end{subfigure}
\caption{Magnitude of the velocity field for two coexisting solutions at $\mu=0.5$ in the bifurcating regime.}
\label{figWALLHUGGING}
\end{figure}

Figure \ref{figWALLHUGGING} illustrates the two qualitatively different and coexisting states in the bifurcation regime at $\mu = 0.5$: the stable asymmetric (\ref{fig:stable_0.5}), and the unstable symmetric (\ref{fig:unstable_0.5}) solutions.
We remark that, for even lower values of $\mu$, the system undergoes a cascade of bifurcations that eventually transition to turbulent and chaotic flow regimes \cite{pintore_2021}. 

\begin{remark}
When dealing with bifurcation problems, it is of utmost importance to be able to reconstruct a specific branch in a neighborhood of the critical parameter value for which uniqueness is lost, i.e.\ to follow a desired branch while varying the value of the parameter \cite{TESIFP}.
In this framework, several \textit{continuation algorithms} \cite{allgower_numerical_1990,KellerLecturesNumericalMethods1988,TESIFP} have been developed in the literature, aiming to generate a sequence of solutions belonging to the same branch.
Here, the continuation algorithm used to retrieve the represented bifurcation diagrams for the steady case of the Navier-Stokes equations is the one presented in~\cite{TESIFP, Khamlich_2022}, which simply consists of choosing the previous solution as the initial guess for Newton's method at each new parameter instance.
\end{remark}

We study the numerical approximation of Equations \eqref{eq:NS1} and \eqref{eq:NS2} with the prescribed initial and boundary conditions for $200$ equispaced values of $\mu\in\mathbb{P}=[0.5, 2]$.

In all simulations, the integration in time stops when the steady state is reached.
Consistent with observations in the steady case, after the bifurcation point $\mu^\ast$, coinciding with the one in the steady case, i.e.\ $\mu^\ast\approx 0.96$, the uniqueness of the solution is lost and three qualitatively different configurations coexist.
For the first one, the steady state solution (meaning, the final configuration reached after time integration has stopped) is fully symmetric and unstable, while the two remaining ones are stable and exhibit a wall-hugging behavior either towards the upper or the lower boundaries.

For all values $\mu\geq\mu^\ast$, the converged steady states align almost perfectly with solutions to the steady problem.
For $\mu<\mu^\ast$, the steady configuration always converges to an asymmetric state, indicating the solver's preference for stable solutions.
Either the final state coincides with the steady-state solution, or it belongs to the opposite stable branch, differing from it only by reflection across the horizontal axis.

Regarding the time evolution of the system before the steady state is reached, the simulations consistently show two different behaviors, before and after the bifurcation point.
For values of $\mu$ preceding the bifurcation point, the evolution of the solution is symmetrical with respect to the horizontal axis, and it stops once the steady state is reached.

\begin{figure}[ht]
\centering

\begin{subfigure}{0.47\linewidth}
\centering
\includegraphics[width=\linewidth, trim=0 50 0 130, clip]{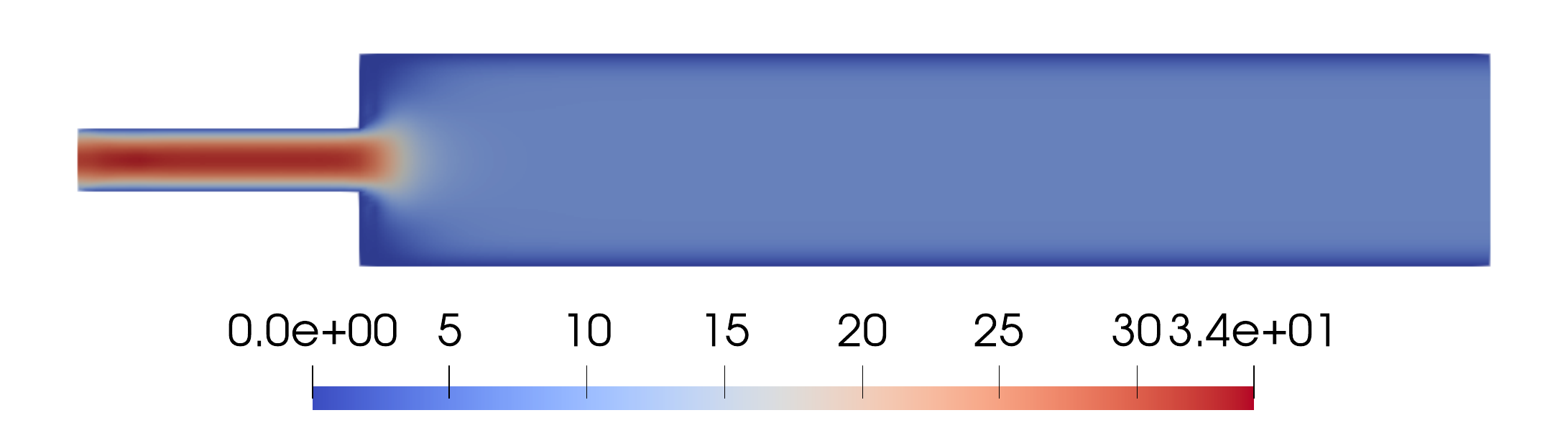}
\caption{$t=0.1$}
\end{subfigure}%
\hfill
\begin{subfigure}{0.47\linewidth}
\centering
\includegraphics[width=\linewidth, trim=0 50 0 130, clip]{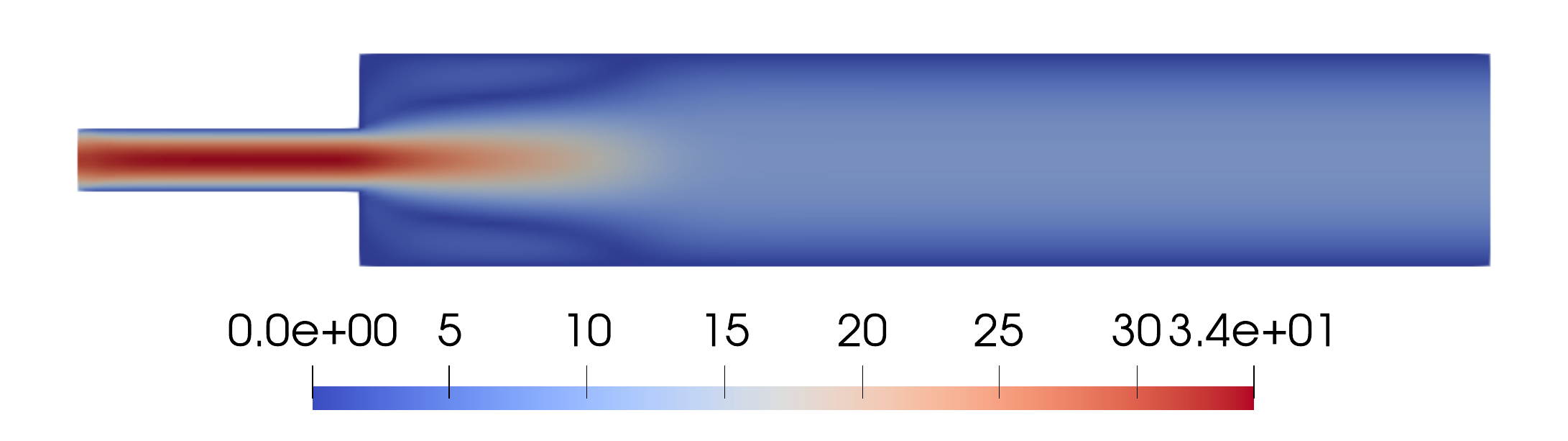}
\caption{$t=1$}
\end{subfigure}

\vspace{0.2cm}

\begin{subfigure}{0.47\linewidth}
\centering
\includegraphics[width=\linewidth, trim=0 50 0 130, clip]{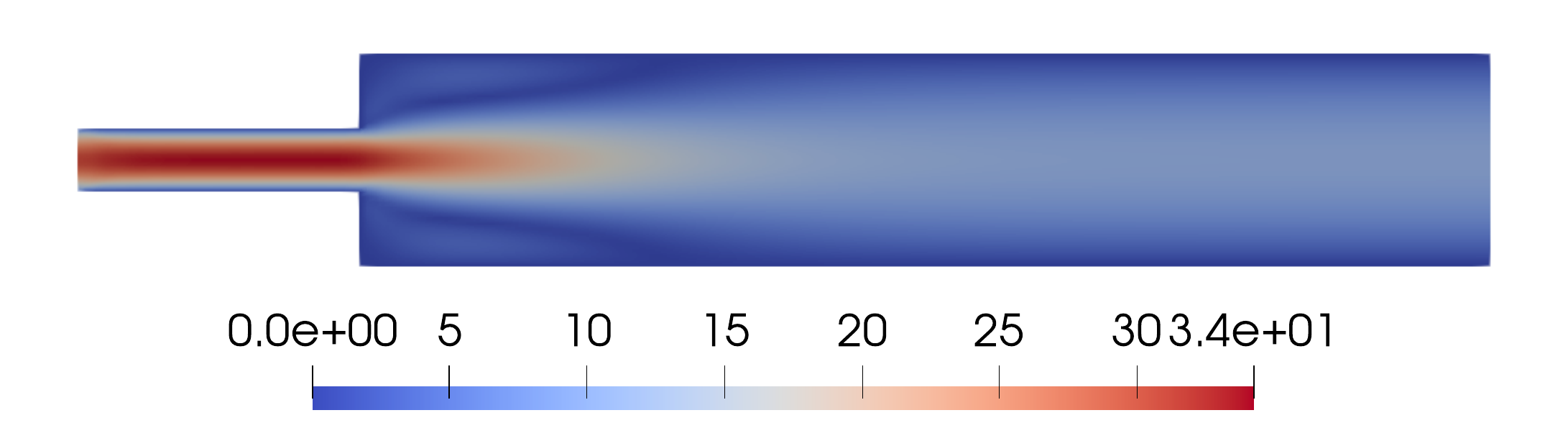}
\caption{$t = T_\text{end}=4.67$}
\end{subfigure}%
\hfill
\begin{subfigure}{0.47\linewidth}
\centering
\includegraphics[width=\linewidth, trim=0 50 0 130, clip]{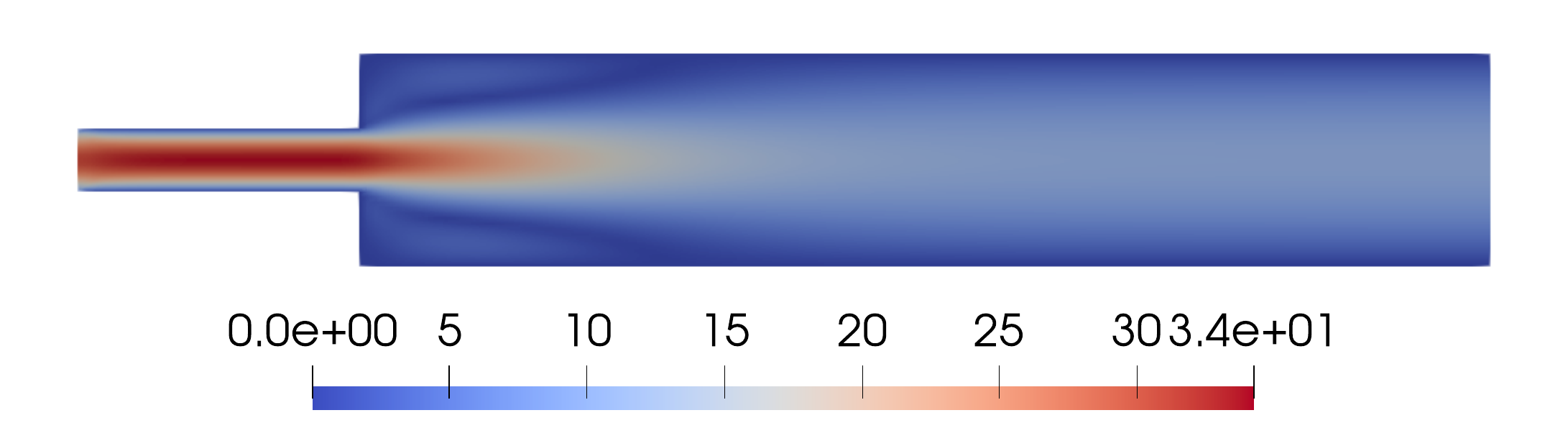}
\caption{Steady solution}
\end{subfigure}

\caption{Evolution of the velocity field's magnitude for $\mu = 1.3$, and comparison with solution of the steady problem.}
    \label{figEvolution1}
\end{figure}

Figure \ref{figEvolution1} illustrates the evolution in time of the flow's velocity for $\mu=1.3$ at several time instances, and the comparison with the solution for the steady problem.

\begin{figure}[ht]
\centering

\begin{subfigure}{0.47\linewidth}
\centering
\includegraphics[width=\linewidth, trim=0 50 0 130, clip]{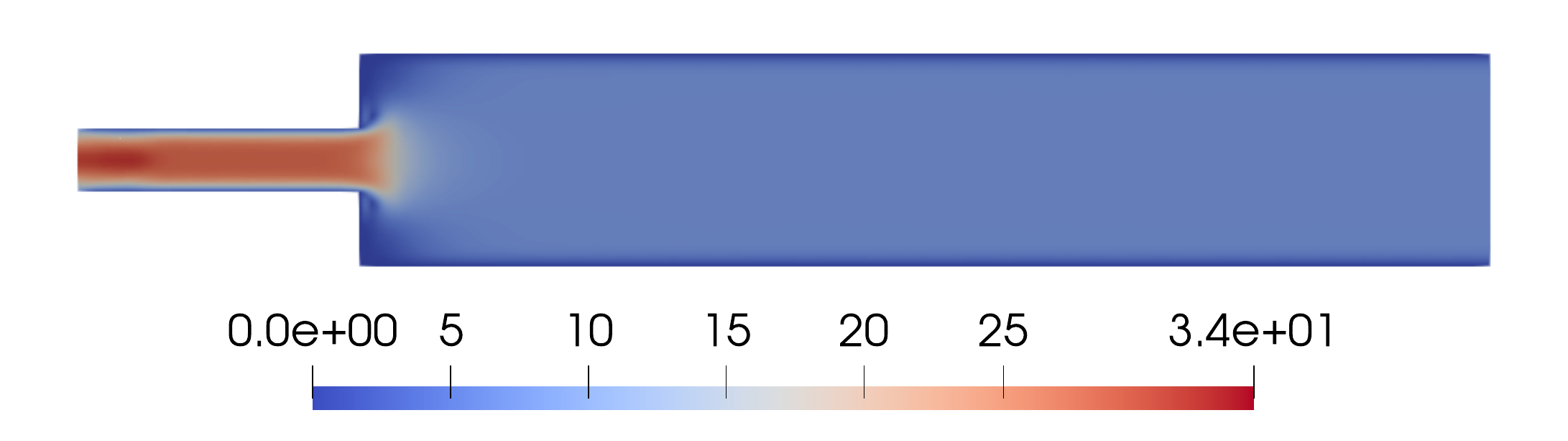}
\caption{$t=0.1$}
\end{subfigure}%
\hfill
\begin{subfigure}{0.47\linewidth}
\centering
\includegraphics[width=\linewidth, trim=0 50 0 130, clip]{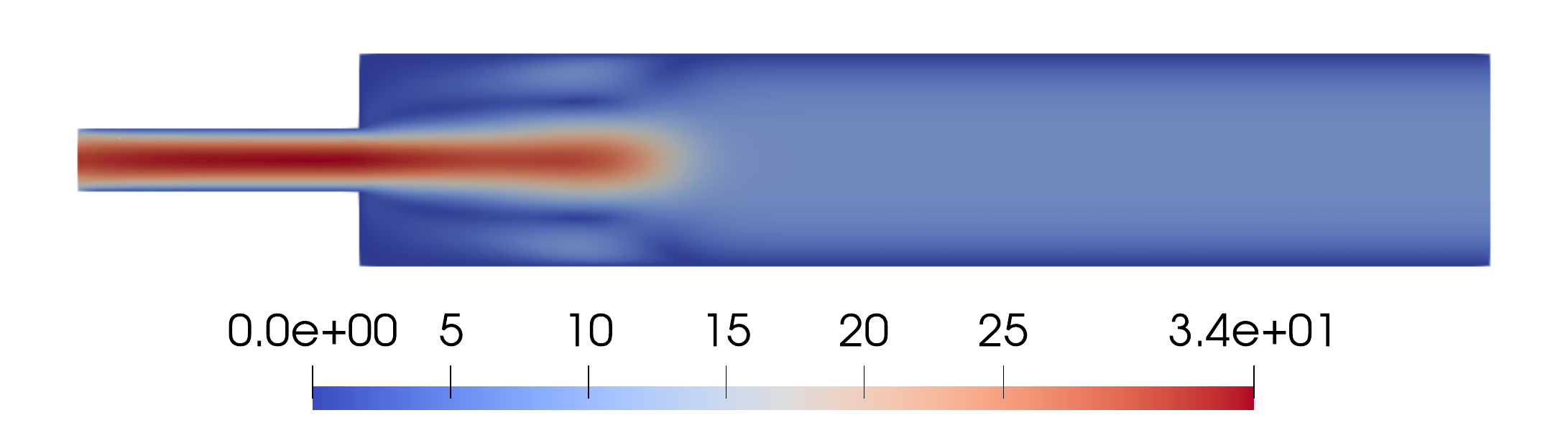}
\caption{$t=1$}
\end{subfigure}

\vspace{0.2cm}

\begin{subfigure}{0.47\linewidth}
\centering
\includegraphics[width=\linewidth, trim=0 50 0 130, clip]{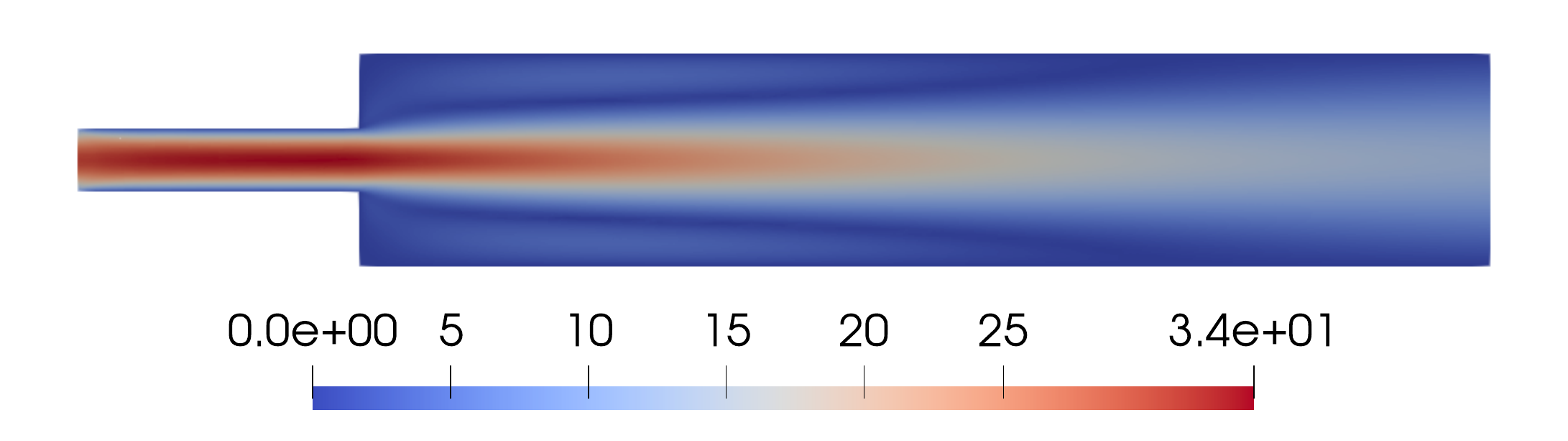}
\caption{$t = 6.31$}
\end{subfigure}%
\hfill
\begin{subfigure}{0.47\linewidth}
\centering
\includegraphics[width=\linewidth, trim=0 50 0 130, clip]{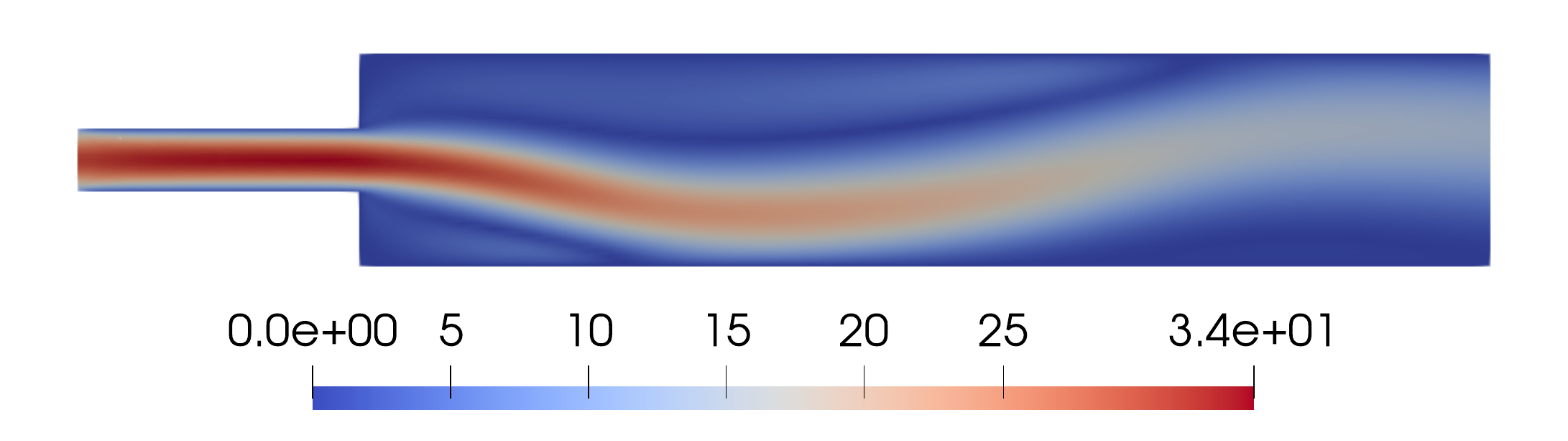}
\caption{$t=T_\text{end}=48.47$}
\end{subfigure}

\vspace{0.2cm}

\begin{subfigure}{0.47\linewidth}
\centering
\includegraphics[width=\linewidth, trim=0 50 0 130, clip]{./images/sec_bifurcations/coanda/steady_symm_mu_0.5.pdf}
\caption{Steady solution wo.\ continuation algorithm}
\end{subfigure}%
\hfill
\begin{subfigure}{0.47\linewidth}
\centering
\includegraphics[width=\linewidth, trim=0 50 0 130, clip]{./images/sec_bifurcations/coanda/steady_asymm_mu_0.5.pdf}
\caption{Steady solution w.\ continuation algorithm}
\end{subfigure}

\caption{Evolution of the velocity field's magnitude for $\mu = 0.5$, and comparison with the solution obtained considering the steady problem.}
    \label{figEvolution2}
\end{figure}

Regarding values of $\mu$ in the bifurcation regime, the evolution of the fluid is symmetric until an almost steady state is reached.
After that moment, the solution begins to bend toward one of the walls, eventually reaching an asymmetric configuration.
Figure~\ref{figEvolution2} displays the evolution of the velocity field for $\mu=0.5$, and the corresponding stable and unstable solutions obtained with or without the exploitation of the continuation algorithm, implemented in~\cite{RozzaBallarinScandurraPichi2024,allgower_numerical_1990}.

\begin{figure}
\begin{subfigure}{0.45\linewidth}
\centering
\includegraphics[width=\linewidth]{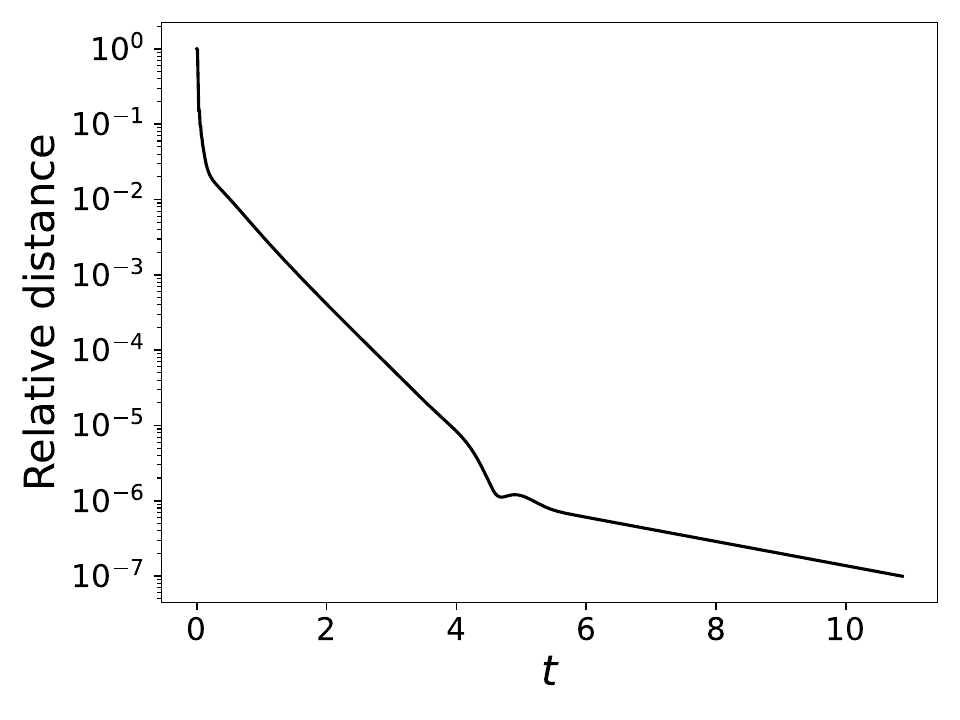}
\caption{$\mu=1.3$}
\end{subfigure}
\hspace{0.6cm}
\begin{subfigure}{0.45\linewidth}
\centering
\includegraphics[width=\linewidth]{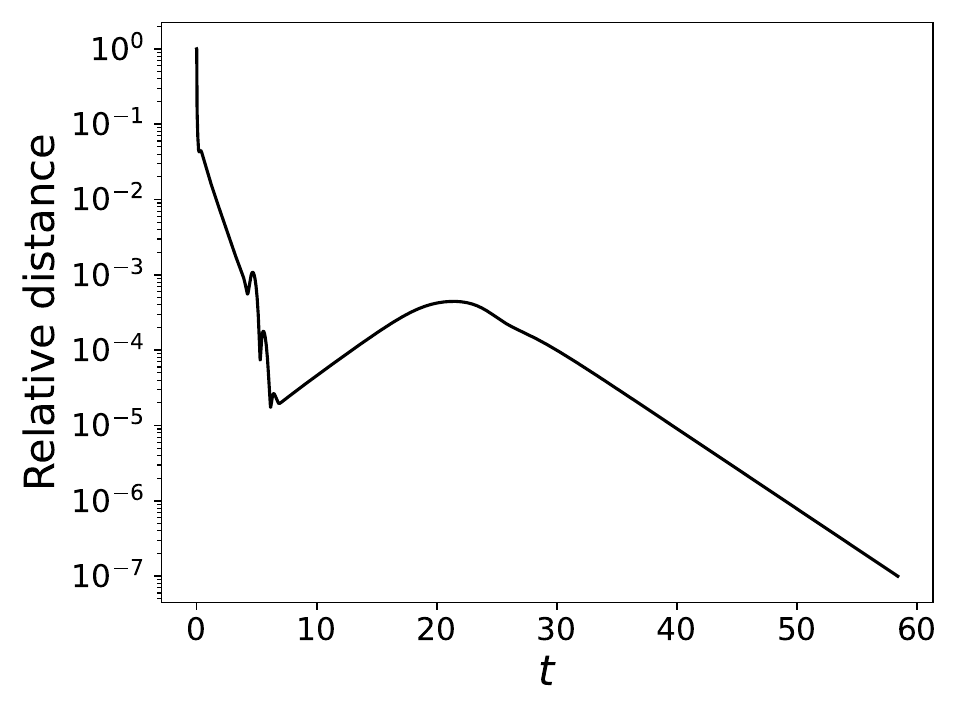}
\caption{$\mu=0.5$}
\end{subfigure}
\caption{Evolution over time of the relative distance between two successive iterations for two values of $\mu$, before and after the bifurcation point.}
    \label{figRELDIST}
\end{figure}

Figure \ref{figRELDIST} shows the evolution of the relative distance between two successive iterations w.r.t.\ time until convergence for $\mu= 1.3$ and $\mu=0.5$.
For values of $\mu\geq\mu^\ast$, the relative distance between two successive iterations decreases almost uniformly until the solver reaches the steady state.
In contrast, for values of $\mu<\mu^\ast$ the relative distance initially decreases, reaching a local minimum in correspondence to the symmetrical configuration.
As the solution begins to bend toward a wall, the relative distance increases before gradually diminishing until it reaches a global minimum corresponding to the steady state.
We also remark that the time required to satisfy the stopping criterion \eqref{eqSTOP} is significantly longer for values of $\mu<\mu^\ast$.

For this reason, it is crucial that the threshold $\tau$ selected to stop the time stepping is sufficiently low to guarantee the full development of the system, and the identification of the branches.
In fact, choosing a threshold that is too high, e.g.\ $10^{-5}$ instead of $10^{-7}$, has a negative impact on the discovery of the bifurcating  states, compromising the proper reconstruction of the bifurcation diagram.

\begin{figure}[ht]
\centering

\begin{subfigure}{0.48\linewidth}
\centering
\includegraphics[width=\linewidth]{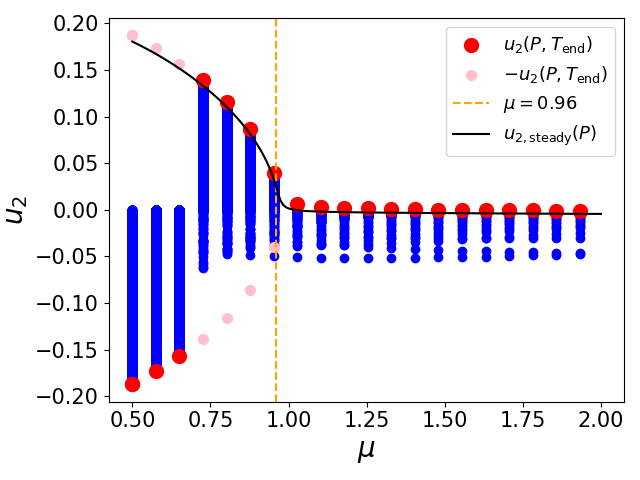}
\caption{Comparison with the steady case}
\label{fig:vert_vel_1}
\end{subfigure}
\hspace{0.2cm}
\begin{subfigure}{0.48\linewidth}
\centering
\includegraphics[width=\linewidth]{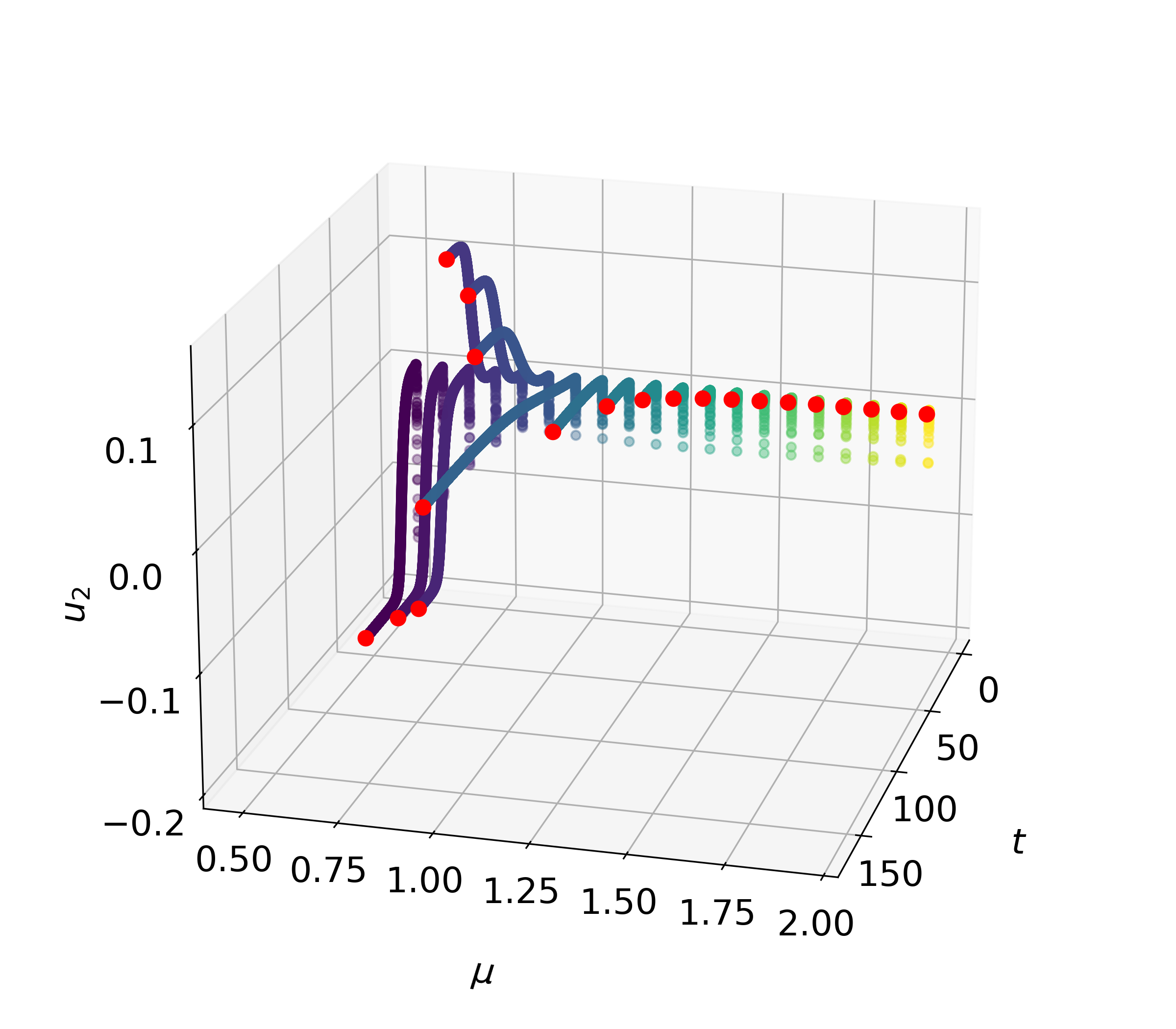}
\caption{3D dynamics}
\label{fig:vert_vel_2}
\end{subfigure}

\caption{Evolution of $u_2$ over time at point $P$ for several values of $\mu$, and comparison with the solution of the steady problem.}
\label{figVERTVEL}
\end{figure}

Having described the behavior of the solutions to Equations \eqref{eq:NS1} and \eqref{eq:NS2}, we are interested in showing the bifurcation diagrams.
This is done by choosing a reasonable quantity of interest $s(\cdot)$ for the parametrized solution at the final time, and plotting its dependence on $\mu$.
A straightforward option is to consider a point on the horizontal axis of symmetry and evaluate the vertical component $u_2$ of the velocity $\bs{u}=(u_1, u_2)$.
To detect asymmetries, we choose the point that lies on the symmetry axis of the geometry with coordinates $P=(19.5, 3.75)$, for which the resulting bifurcation diagram is plotted in Figure \ref{figVERTVEL}.
In Figure \ref{fig:vert_vel_1}, we observe that once a steady state is reached, the solution's behavior almost perfectly overlaps with the one of the solution to the steady equation.
Furthermore, for each value of $\mu$ in the bifurcation regime, the solver converges to a stable asymmetric solution displaying wall-hugging behavior; thus, for consistency, we plotted the evolution of both stable branches.
Inspecting Figure \ref{fig:vert_vel_2}, we notice that the number of time step needed to reach the steady state can indicate a-posteriori the location of the bifurcation point.

\begin{figure}[ht]
\centering
\begin{subfigure}{0.48\linewidth}
\centering
\includegraphics[width=\linewidth]{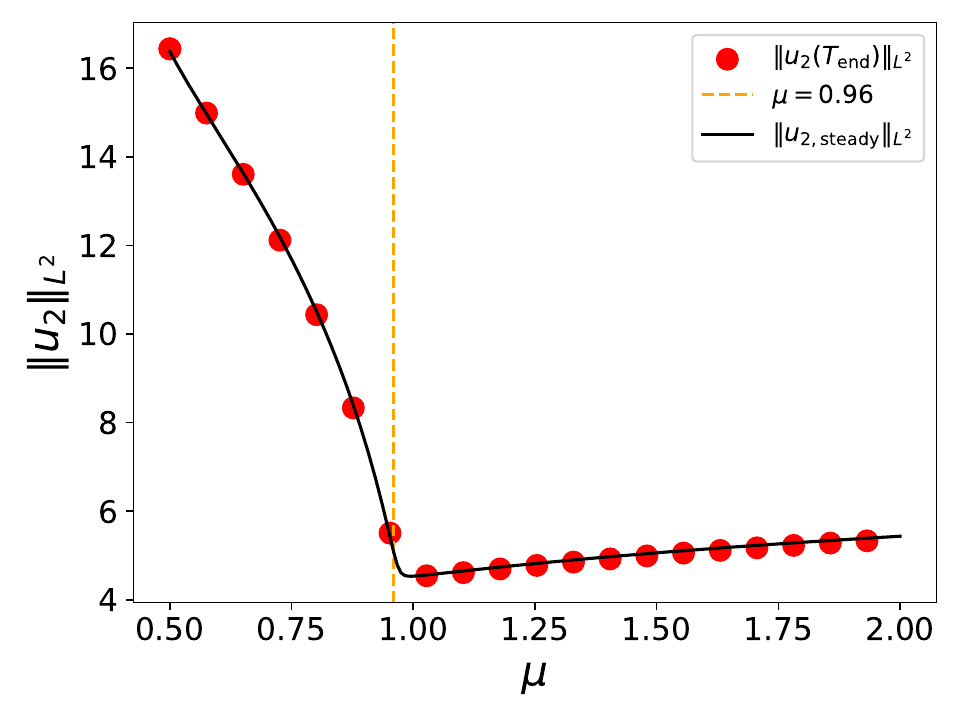}
\caption{Comparison with the steady case}
\label{fig:norm_1}
\end{subfigure}
\vspace{0.6cm}
\begin{subfigure}{0.48\linewidth}
\centering
\includegraphics[width=\linewidth]{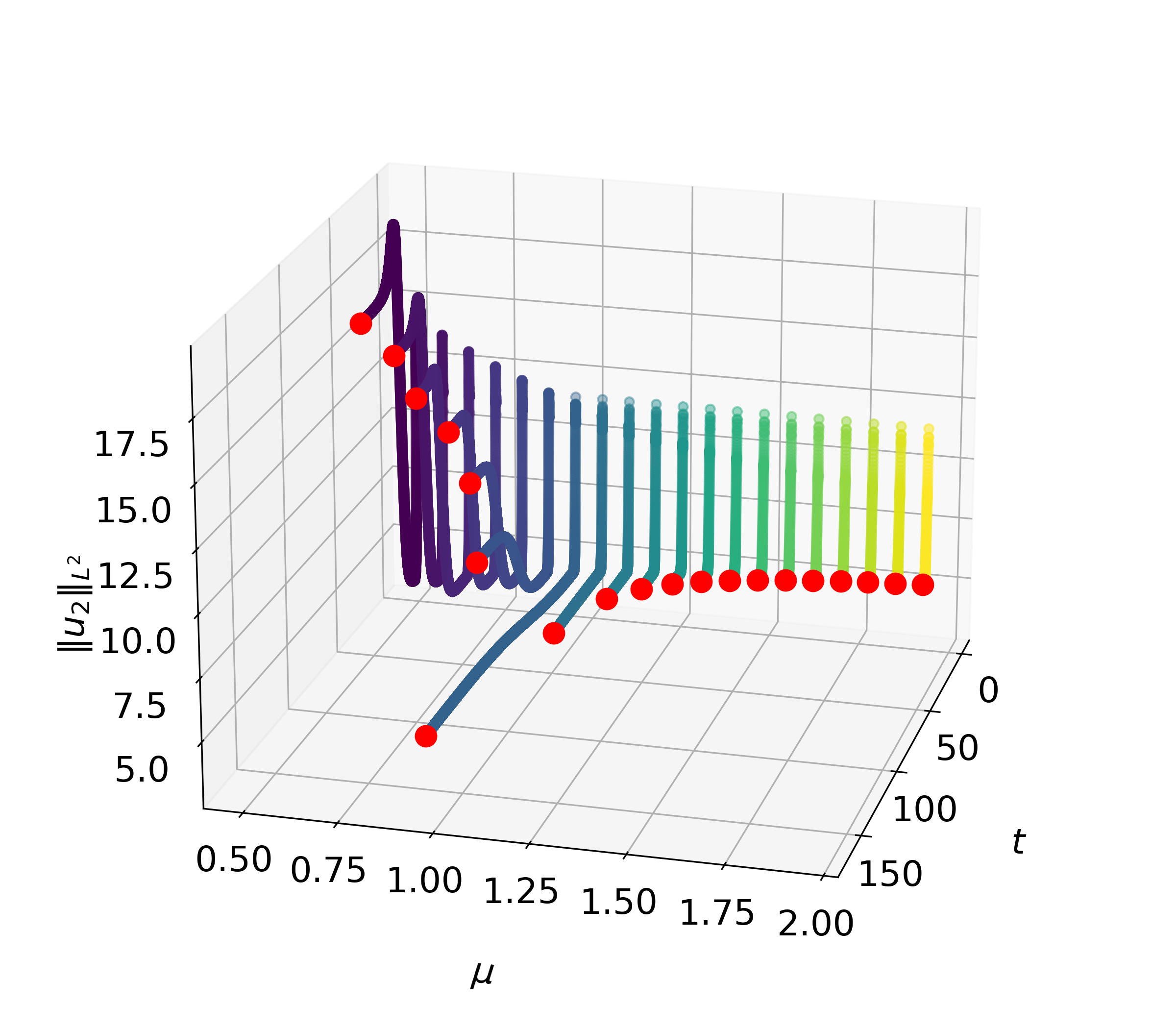}
\caption{3D dynamics}
\label{fig:norm_2}
\end{subfigure}

\caption{Evolution over time of the $L^2$-norm of the vertical component of the velocity field $u_2$ for several values of $\mu$, and comparison with the solution of the steady problem.}
\label{figNORMy}
\end{figure}

Exploiting the $L^2-$norm of the vertical velocity $u_2$ is very informative and allows reconstructing the bifurcation diagram, as can be seen in Figure \ref{fig:norm_1}.
In addition, Figure \ref{fig:norm_2} displays the number of iterations needed to reach the steady state. 
This confirms that the main cause of the bifurcation is the symmetry breaking that occurs for the vertical component of the velocity.

\begin{remark}
Varying the parameter $\theta$ had minimal impact on the overall behavior, though $\theta = 1$ resulted in longer transient phases at low $\mu$ values before wall-hugging behavior emerged.
To obtain stable solutions directly and avoid the symmetrical transient phase, we explored multiple approaches aimed at choosing a suitable initial guess for the nonlinear solver at each time step. These included: (i) using asymmetrical inflow conditions to generate a suitable initial guess as done in \cite{Hess_2019}, and (ii) using solutions from the continuation algorithm or their weighted combinations with previous time steps as initial guesses for the nonlinear solver.
Neither approach altered the evolution of the solution\footnote{
Additional numerical experiments can be found at  
\href{https://github.com/lorenzotomada/NS_continuation}{https://github.com/lorenzotomada/NS\_continuation}.}.
\end{remark}

\subsection{Hopf bifurcation}\label{subsec:hopf_problem}
\begin{figure}[ht]
    \centering
    \includegraphics[width=1\textwidth]{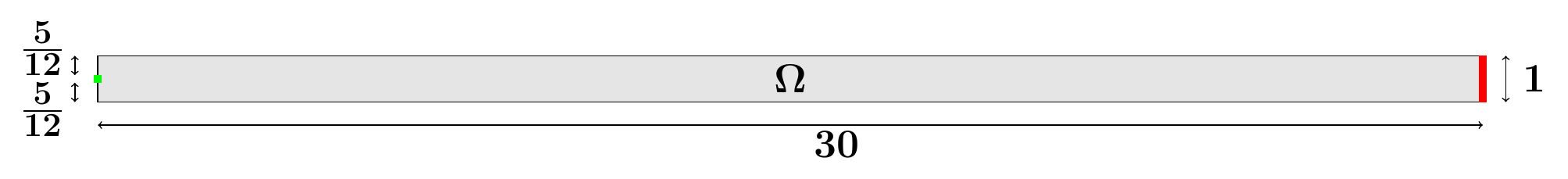}
    \caption{Problem geometry showing inlet $\Gamma_\text{in}$ (green), outlet $\Gamma_\text{out}$ (red), and wall boundary $\Gamma_0$ (black).}
    \label{figQUAINI}
\end{figure}

Our second case study examines fluid flow in a sudden-expansion channel. 
The domain $\Omega = (0,1)\times(0,30)$ is shown in Figure \ref{figQUAINI}, with boundary conditions:
\begin{equation}\label{eq:bc_hopf}
    \begin{dcases}
    \bs{u}(\bs{x}, 0) = \bs{0}, &\text{in } \Omega, \\
    \bs{u}(\bs{x}, t) = [U_{\text{max}}(1 - 4(y - \frac{1}{2})^2), 0], &\text{on } \Gamma_{\text{in}}\times(0,T), \\
    \bs{u}(\bs{x}, t) = \bs{0}, &\text{on } \Gamma_0\times(0,T), \\
    (-\mu\nabla\bs{u} + p\mathbf{I})\bs{n} = \bs{0}, &\text{on } \Gamma_{\text{out}}\times(0,T),
    \end{dcases}
\end{equation}
where $\Gamma_{\text{in}}=\{0\}\times [\frac{5}{12},\,\frac{7}{12}]$, $\Gamma_{\text{out}} = \{30\}\times[0,\, 1]$, and $\Gamma_0 = \partial\Omega\smallsetminus(\Gamma_{\text{in}}\cup\Gamma_{\text{out}})$.

Following \cite{QUAINI}, we define the Reynolds number $Re$ as our bifurcation parameter:
\begin{equation}\label{eq:RE_QUAINI}
Re = \frac{2U\rho L}{\mu},
\end{equation}
where $U$ represents the average channel velocity, $\rho$ denotes fluid density, and $L=\frac{1}{6}$ is the characteristic length. Given the maximum velocity $U_{\text{max}}=1$ and assuming a fully developed parabolic profile, we have $U=\frac{2}{3}U_\text{max}$, yielding:
\begin{equation}\label{eq:RE_QUAINI2}
    Re = \frac{\rho U}{3\mu} = \frac{2}{9}\frac{\rho U_\text{max}}{\mu}.
\end{equation}

\begin{figure}[ht]
\centering

\begin{subfigure}{0.5\linewidth}
\centering
\includegraphics[width=\linewidth, trim=350 180 350 230, clip]{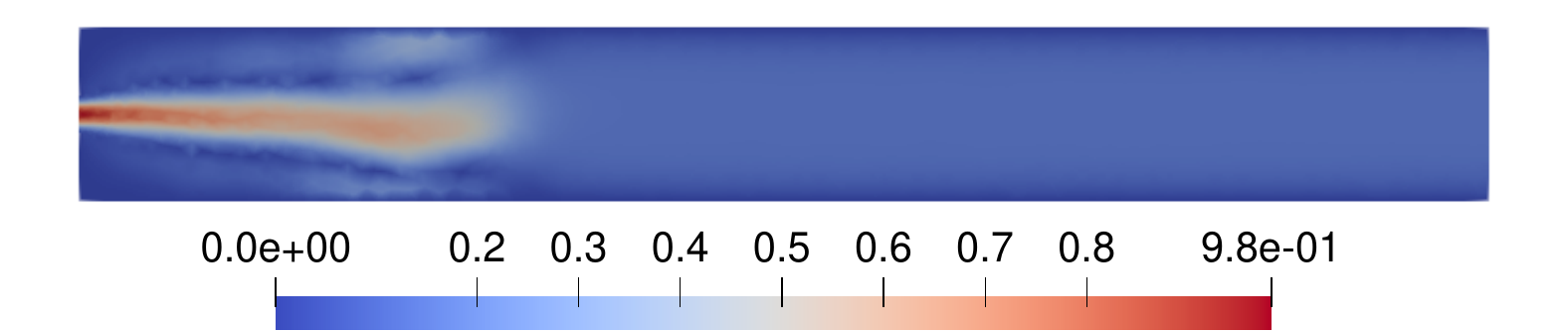}
\caption{$t=10$}
\end{subfigure}%
\hfill
\begin{subfigure}{0.5\linewidth}
\centering
\includegraphics[width=\linewidth, trim=350 180 350 230, clip]{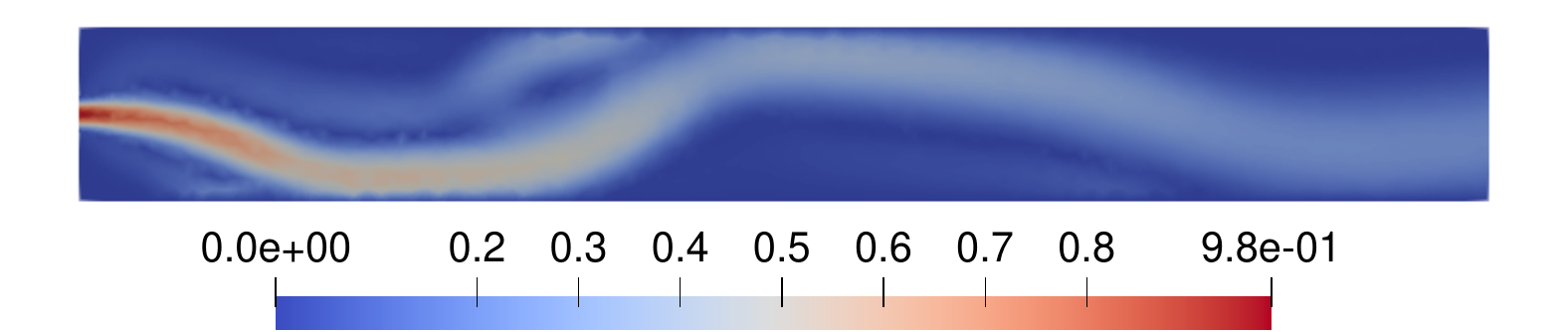}
\caption{$t=100$}
\end{subfigure}

\begin{subfigure}{0.5\linewidth}
\centering
\includegraphics[width=\linewidth, trim=350 180 350 230, clip]{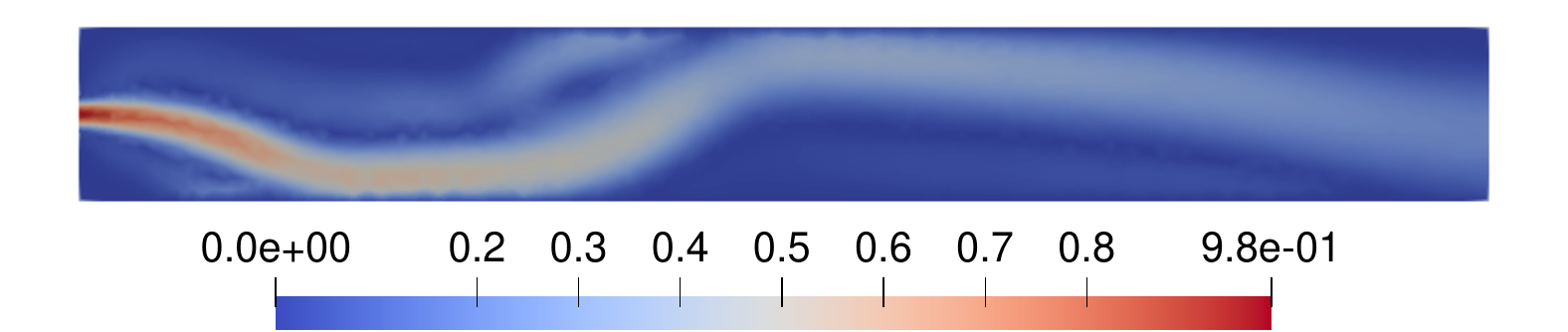}
\caption{$t=200$}
\end{subfigure}%
\hfill
\begin{subfigure}{0.5\linewidth}
\centering
\includegraphics[width=\linewidth, trim=350 180 350 230, clip]{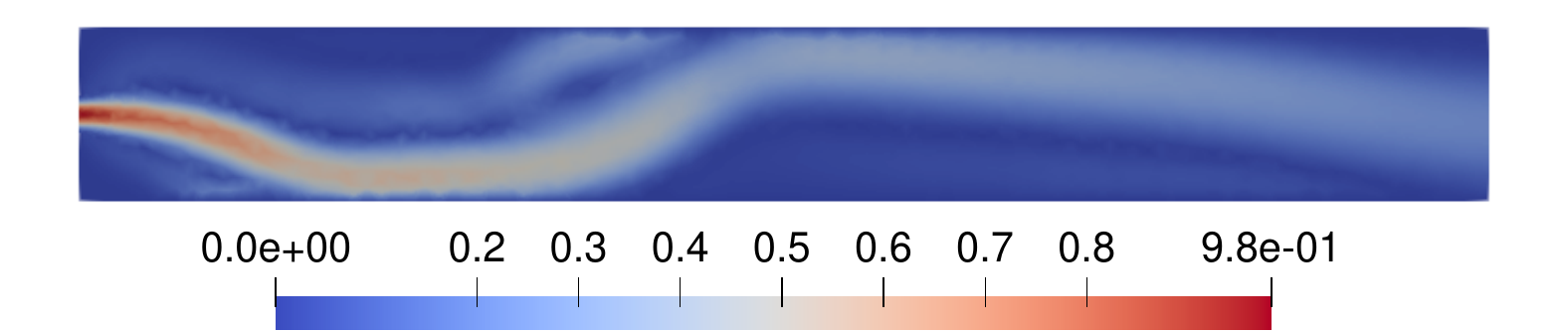}
\caption{$t=300$}
\end{subfigure}

\caption{Evolution of the velocity field's magnitude for $Re=346$ (before bifurcation), showing steady behavior. The flow quickly develops into a steady-state stable asymmetric pattern, with the jet attaching to one wall due to the Coandă effect.}
\label{figEvolutionHopf346}
\end{figure}

\begin{figure}[ht]
\centering

\begin{subfigure}{0.5\linewidth}
\centering
\includegraphics[width=\linewidth, trim=350 180 350 230, clip]{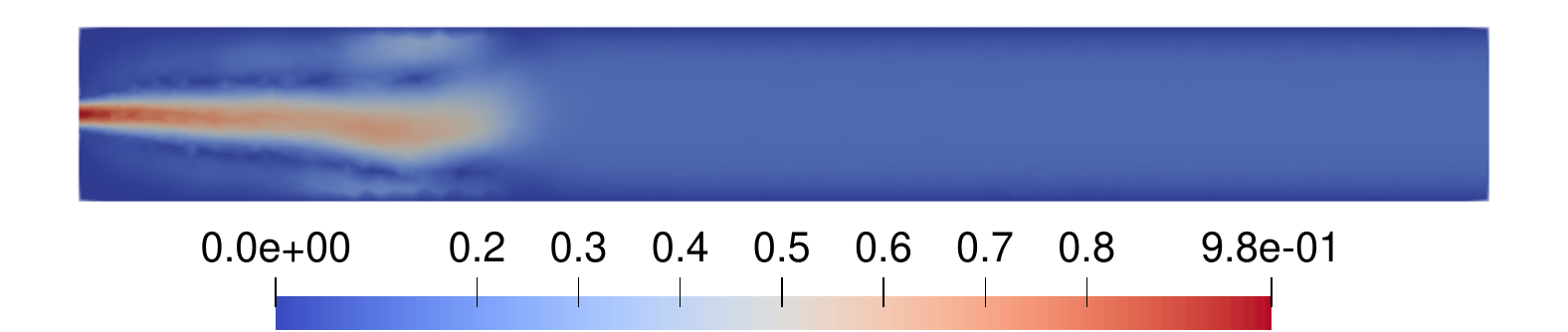}
\caption{$t=10$}
\end{subfigure}%
\hfill
\begin{subfigure}{0.5\linewidth}
\centering
\includegraphics[width=\linewidth, trim=350 180 350 230, clip]{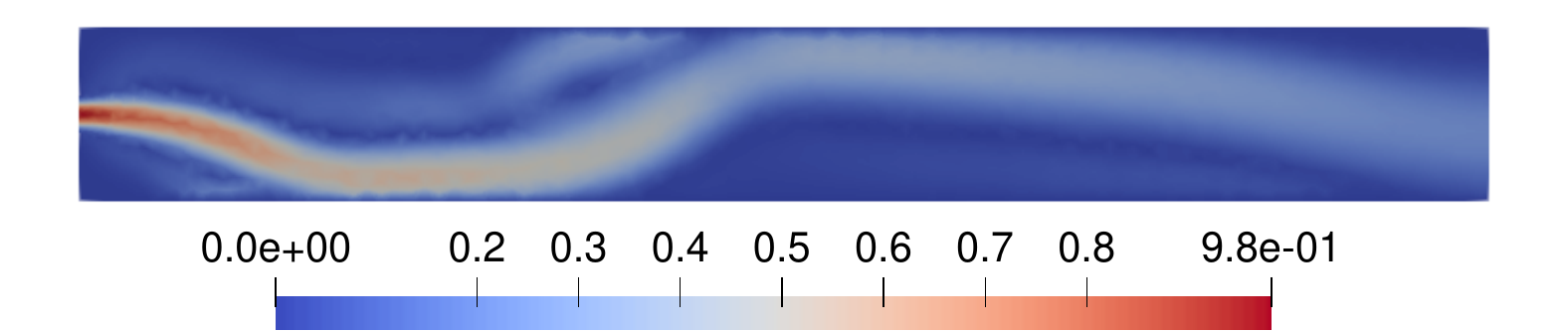}
\caption{$t=200$}
\end{subfigure}

\vspace{0.2cm}

\begin{subfigure}{0.5\linewidth}
\centering
\includegraphics[width=\linewidth, trim=350 180 350 230, clip]{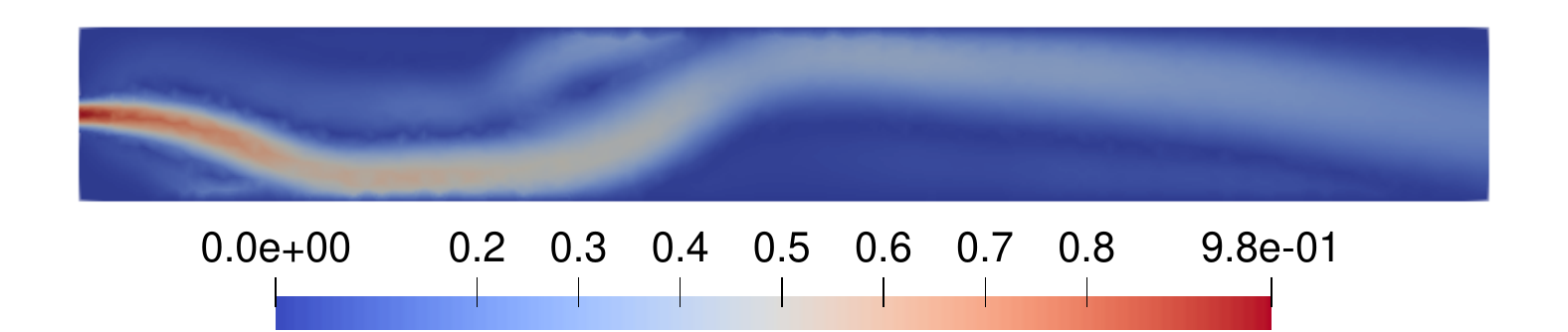}
\caption{$t=400$}
\end{subfigure}%
\hfill
\begin{subfigure}{0.5\linewidth}
\centering
\includegraphics[width=\linewidth, trim=350 180 350 230, clip]{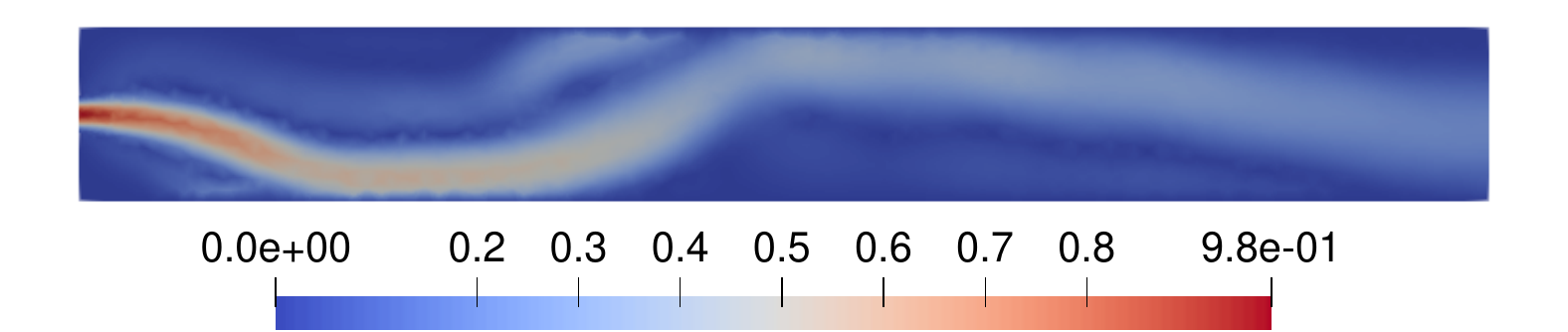}
\caption{$t=600$}
\end{subfigure}

\begin{subfigure}{0.5\linewidth}
\centering
\includegraphics[width=\linewidth, trim=350 180 350 230, clip]{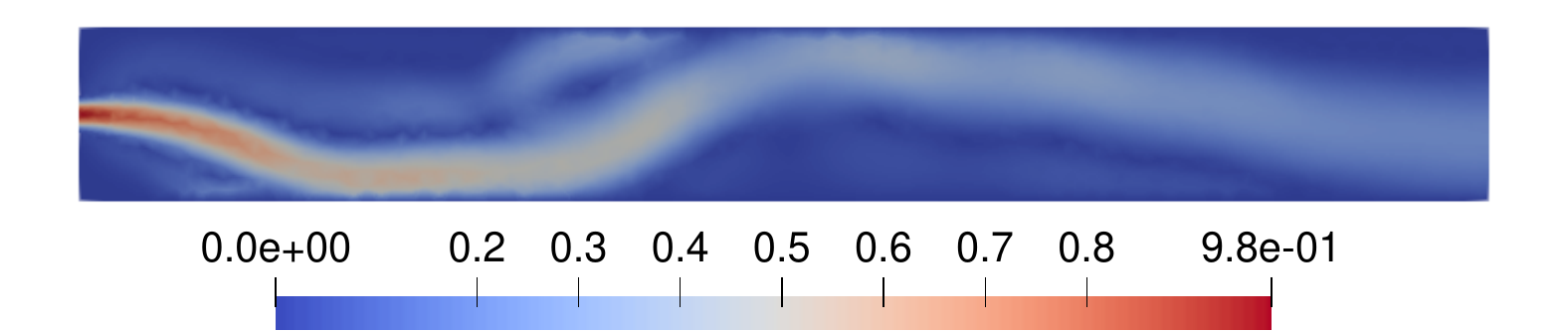}
\caption{$t=800$}
\end{subfigure}%
\hfill
\begin{subfigure}{0.5\linewidth}
\centering
\includegraphics[width=\linewidth, trim=350 180 350 230, clip]{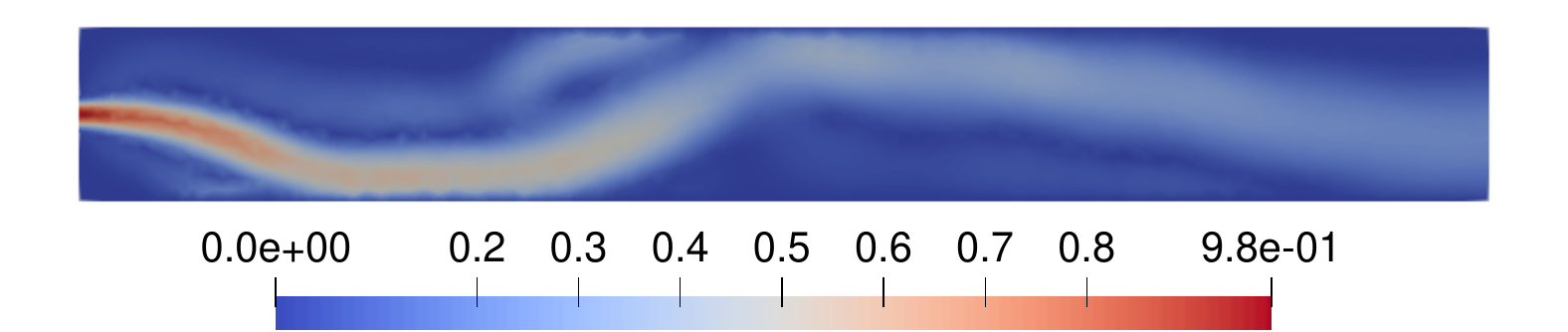}
\caption{$t=900$}
\end{subfigure}

\caption{Evolution of the velocity field's magnitude for $Re=356$ (after bifurcation), showing the emergence of periodic behavior. The initially asymmetric flow develops instabilities that lead to periodic vortex shedding, reaching the characteristic limit cycle behavior.}
\label{figEvolutionHopf356}
\end{figure}
For the numerical simulations, we focus on the parameter range $Re\in\mathbb{P}=[346, 356]$, which encompasses the Hopf bifurcation point, and we employ the \gls{fv} method implemented in OpenFOAM \cite{Weller_1998}. The pressure-velocity coupling is handled through the PIMPLE algorithm \cite{Greenshields_2022}, while spatial discretization utilizes second-order central differences for diffusion terms and second-order upwind schemes for convection terms. Our mesh consists of $\num{12734}$ cells, which is intentionally coarser than those used in \cite{QUAINI}. Time integration employs a stepsize of $dt=10^{-2}$, though solution data is recorded at 1-second intervals to optimize memory usage, resulting in a temporal discretization $I_k=[k\mid k\in\mathbb{N},\, k\leq 2500]$.

\subsubsection{Bifurcation Analysis}

The system exhibits a rich progression of dynamical behaviors as the Reynolds number increases. At low Reynolds numbers, the flow transitions from symmetric to asymmetric solutions through the Coandă effect, similar to the pitchfork bifurcation case. As $Re$ increases further, the system develops four distinct stable asymmetric solutions. The most interesting transition occurs at higher Reynolds numbers, where these asymmetric solutions lose stability through a Hopf bifurcation, giving rise to periodic solutions \cite{QUAINI}. Beyond this critical point, the elongated jet becomes unstable, fragmenting into small vortices that propagate downstream, establishing a time-periodic flow pattern. This behavior is clearly visible in the velocity field visualizations shown in Figures \ref{figEvolutionHopf346} and \ref{figEvolutionHopf356}, which contrast the steady-state behavior before bifurcation with the periodic dynamics after bifurcation.
For visualization purposes, only the initial portion of the channel (with $x\in[0,8]$) is shown, as the most relevant dynamics occur in this region.

\begin{figure}[ht]
	\centering
    \includegraphics[width=1.0\textwidth]{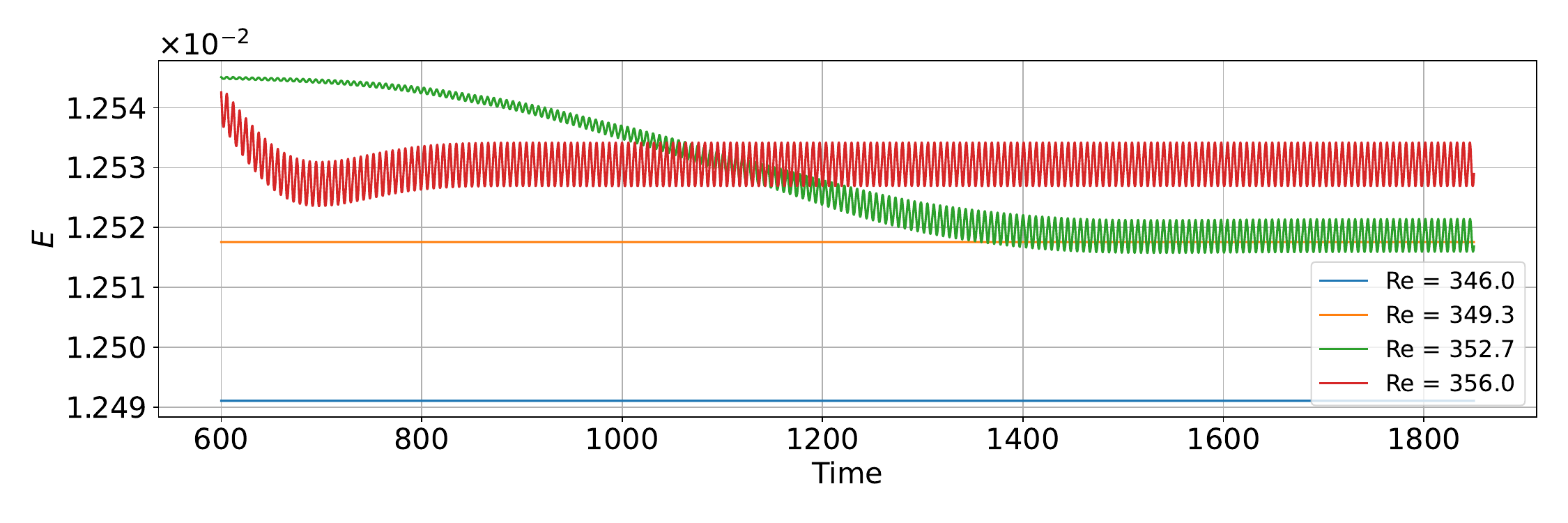}
	\caption{Evolution of the kinetic energy $E(t)$ for selected values of $Re$, demonstrating the transition from steady to periodic behavior.}
	\label{figENERGY1}
\end{figure}

To systematically investigate this behavior and the evolution in time of the system, we consider $40$ equispaced values of $Re$ in  $\mathbb{P}=[346, 356]$, and 
we monitor the total kinetic energy of the system, defined as:
\begin{equation*}
E(t) = \frac{1}{2}\int_\Omega|\bs{u}_h(t)|^2d\Omega,
\end{equation*}
which serves as an effective indicator of the system's dynamics. Figure \ref{figENERGY1} presents the evolution of $E(t)$ over the interval $[600, 1850]$ for four carefully selected Reynolds numbers—two before and two after the bifurcation point—illustrating the transition in system behavior.

\begin{figure}[ht]
	\centering
    \includegraphics[width=0.48\textwidth]{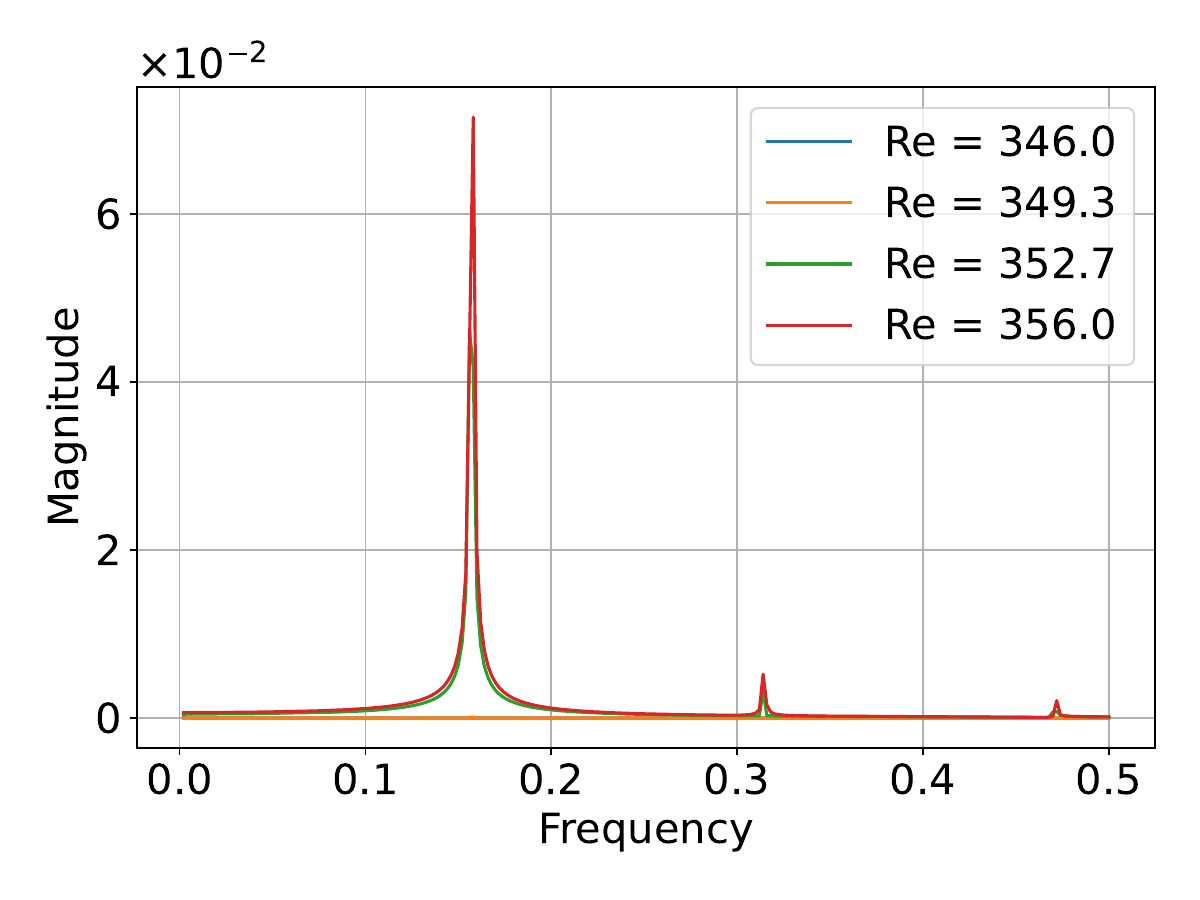}
	\caption{Power spectral density analysis of $E(t)$ revealing the frequency content of the flow for selected Reynolds numbers.}
	\label{figFFT_energy}
\end{figure}

The periodic nature of the post-bifurcation solutions can be rigorously established through Fourier analysis of the kinetic energy signal $E(t)$. Figure \ref{figFFT_energy} displays the power spectral density of $E(t)$ for the same four Reynolds numbers, computed from time series data comprising $\num{50000}$ points\footnote{Unlike the full-order snapshots, the energy values were saved every $10^{-2}$ seconds, allowing for a more complete investigation of the bifurcating phenomenon.} over $[\num{2000}, \num{2500}]$. The emergence of distinct frequency peaks after the bifurcation point clearly indicates the onset of periodic behavior.

\begin{figure}[ht]
    \centering 
    \includegraphics[width=.55\textwidth]{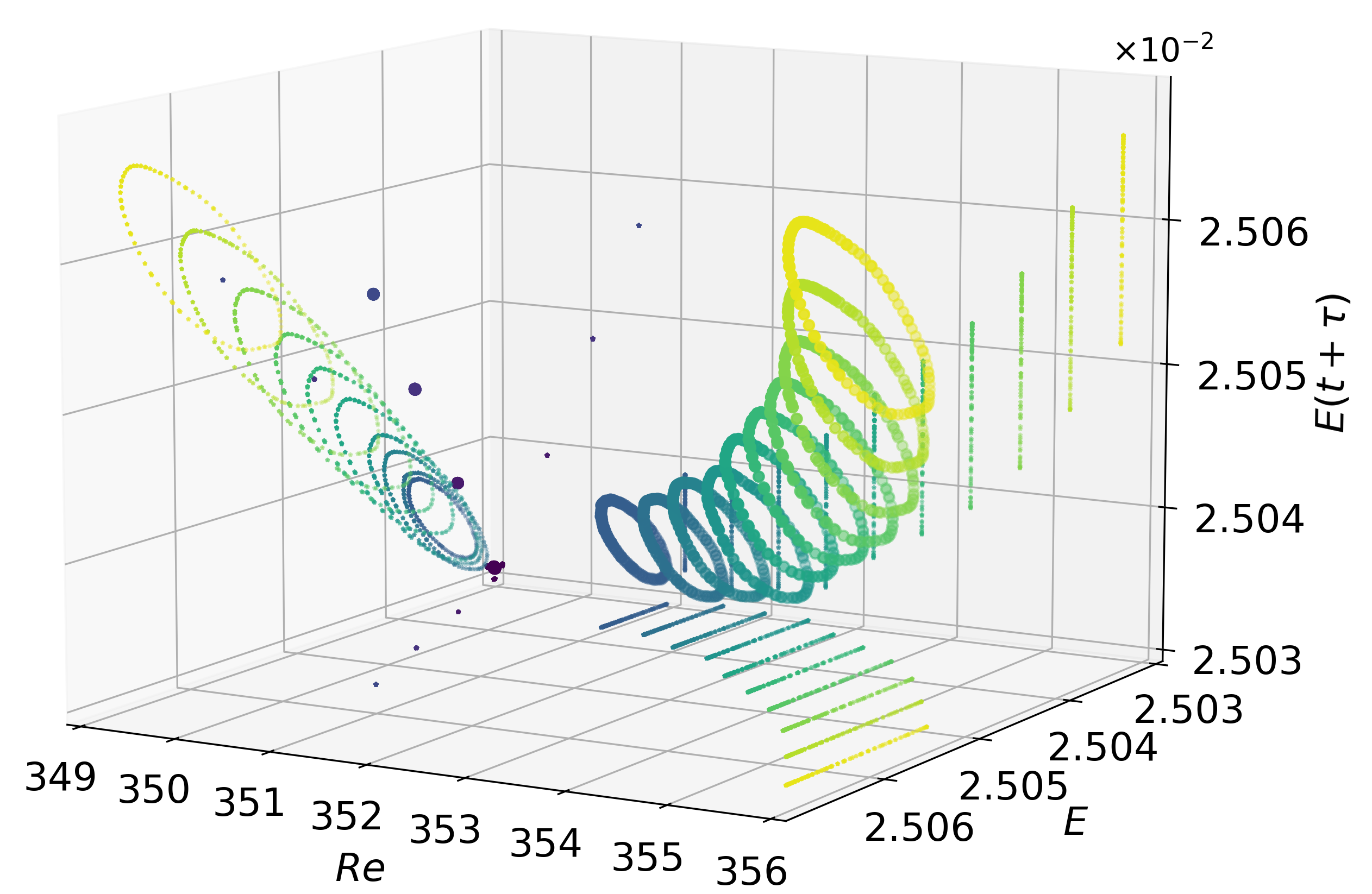}
    \caption{3D phase space representations of the system dynamics across different Reynolds numbers.}
    \label{fig:energy}
\end{figure}

Moreover, the periodic solutions can also be represented in the three-dimensional phase plane.
Figure \ref{fig:energy} presents a plot of $E(t)$ against $E(t+\tau)$ across multiple Reynolds numbers, with $\tau=7$ and starting from $t=1800$.
This phase portrait clearly demonstrates the transition from fixed points to limit cycles, characteristic of a Hopf bifurcation.

To quantitatively characterize the bifurcation, we introduce a scalar measure $A: \mathbb{P}\to\R^+$ representing the amplitude of energy oscillations:
\begin{equation}\label{eqAMPLITUDE}
A(Re) = \max_{t\in I_k} E(t) - \min_{t\in I_k} E(t),
\end{equation}
where the dependence of $E$ on $Re$ is implicit.
The resulting bifurcation diagram, presented in Figure \ref{figAMPLITUDE_PARAM} of Section \ref{sec:numerical}, shows the supercritical nature of the Hopf bifurcation.

\begin{remark}\label{rmkBIF_QUAINI}
The location of the bifurcation point exhibits significant mesh sensitivity, as noted in \cite{QUAINI}. Our computed critical value $Re^\ast\approx 350.6$ is lower than the range $[413, 476]$ reported in the original study, which employed finer spatial discretizations. While this discrepancy stems from our intentionally coarser mesh, it does not impact the fundamental bifurcation structure or the effectiveness of our model reduction strategies, which are the primary focus of this investigation.
\end{remark}

\section{Reduced order modeling: SINDy, AE and nested POD}\label{sec:ROM}
The need to mitigate the computational burden when solving parameterized PDEs for multiple parameter values, coupled with the necessity for real-time responses, has driven the development of \gls{rom} techniques.

A typical approach to build reduced models starts by performing expensive computations for some representative states of the system, i.e.\ high-fidelity solutions of the problem, or \textit{snapshots},  upon which a simplified model is built. Given a differential problem with solution $\bs{u}$ and its numerical approximation $\bs{u}_h$ obtained by a \gls{fom}, the \gls{rom} should provide an approximation $\bs{u}_N$ such that $\|\bs{u}_N-\bs{u}\|\approx\|\bs{u}_h-\bs{u}\|$ for a suitably chosen norm.

This section introduces the tools used for developing the \gls{rom} strategies presented in this work, focusing on the combination of SINDy for dynamics identification, autoencoders for nonlinear dimensionality reduction, and nested POD for efficient basis computation.

\subsection{Sparse Identification of Nonlinear Dynamics}\label{subsec:sindy}
This subsection is dedicated to the description of the data-driven method \gls{sindy} as applied to parameterized PDEs, and to some of the critical issues that might arise in the aforementioned situation \cite{Brunton_2016,Champion_2019,SINDyCP,Conti_2023}.
We recall that \gls{sindy} leverages sparse regression with the aim of explicitly retrieving the equations that describe the dynamical system under consideration.
In order to do so, it takes as input some snapshots of the dynamics and a library of candidate functions for its description; the equations are then obtained by numerically solving a sparse regression problem.
Assume that we have a set of snapshots describing the behavior of a parameterized system stored in the snapshot matrix:
\begin{equation}\label{eqSINDyParam} 
\bs{X}=
\begin{bmatrix}
\bs{x}(t_1; \bs{\mu}_1)^T \\
\bs{x}(t_2; \bs{\mu}_1)^T \\
\vdots \\
\bs{x}(t_{N_t}; \bs{\mu}_1)^T \\
\bs{x}(t_1; \bs{\mu}_2)^T \\
\vdots \\
\bs{x}(t_{N_t}; \bs{\mu}_{N_{\mu}})^T
\end{bmatrix}
=
\begin{bmatrix}
x_1(t_1;\bs{\mu}_1) & x_2(t_1;\bs{\mu}_1) & \cdots & x_N(t_1;\bs{\mu}_1) \\
x_1(t_2;\bs{\mu}_1) & x_2(t_2;\bs{\mu}_1) & \cdots & x_N(t_2;\bs{\mu}_1) \\
\vdots & \vdots & \ddots & \vdots\\
x_1(t_{N_t};\bs{\mu}_1) & x_2(t_{N_t};\bs{\mu}_1) & \cdots & x_N(t_{N_t};\bs{\mu}_1) \\
x_1(t_1;\bs{\mu}_2) & x_2(t_1;\bs{\mu}_2) & \cdots & x_N(t_1;\bs{\mu}_2) \\
\vdots & \vdots & \ddots & \vdots\\
x_1(t_{N_t};\bs{\mu}_{N_\mu}) & x_2(t_{N_t};\bs{\mu}_{N_\mu}) & \cdots & x_N(t_{N_t};\bs{\mu}_{N_\mu}) 
\end{bmatrix}
\in\mathbb{R}^{N_tN_{\mu}\times N_h},
\end{equation}
where $N_{\mu}$ is the number of samples of the parameter(s) on which the system depends, and $N_t$ is the number of time instances.
The vector $\bs{x}(t_i;\bs{\mu}_j) = [x_1(t_i;\bs{\mu}_j),\, x_2(t_i;\bs{\mu}_j),\, \cdots,\, x_N(t_i;\bs{\mu}_j)]^T\in\R^{N_h}$, for all $i\in\{1,\dots, N_t\}$ and $j\in\{1,\dots, N_\mu\}$, represents the state of the system at time $t_i$, with $N_h$ being the spatial dimension of the system, while $\bs{\mu}_j\in\R^p$ is a vector of $p$ parameters and/or forcing terms.

Our goal is to learn the dynamics from the data, i.e.\ find a suitable function $\bs{f}$ allowing us to write explicitly a parameterized dynamical system of the form 
\begin{equation}\label{eqParamDyn}
    \begin{cases}
    \dot{\bs{x}}(t; \bs{\mu}) = \bs{f}(t, \bs{x}(t; \bs{\mu}); \bs{\mu})\\[1mm]
    \bs{x}(t_0; \bs{\mu}) = \bs{x}_{0}
    \end{cases}
\end{equation}
where $I\subseteq\R^+$ is some interval, $\bs{f}(\cdot\,;\bs{\mu}):I\times\Omega\to\R^{N_h}$, $\bs{x}(\cdot\,; \bs{\mu}):I\to\Omega$, and $(t_0, \bs{x}_0)\in I\times\Omega$.
The identified model $\bs{f}(t, \bs{x}(t;\bs{\mu});\bs{\mu})$ provides an explicit  formulation in terms of the parameter $\bs{\mu}$, allowing for a comprehensive tool to investigate the whole parametric range.

In order to achieve this, we assume that $\bs{f}$ can be expressed as a linear combination of terms from a library of $r$ candidate functions\footnote{We assume that all candidate functions of both the state variables and the parameters are continuous. While one could consider a library including discontinuous functions, these must provide pointwise evaluations.} $\bs{\Theta}(\bs{x};\,\bs{\mu})\in(C^0[\R^{N_h+p},\R])^r$, and that only a few terms are active i.e.\ the representation is sparse. Hence, we are looking for coefficients $\{\bs{\xi}_k\}_{k=1}^{N_h}$, with each $\bs{\xi}_k\in\R^r$, allowing us to write
\begin{equation}\label{eq:symbolic_SINDy}
    \dot{x}_k(t;\bs{\mu}) = f_k(t, \bs{x}(t;\bs{\mu});\bs{\mu}) = \bs{\Theta}\bigl(\bs{x}(t;\bs{\mu});\bs{\mu}\bigr)\,\bs{\xi}_k, \quad k = 1, \dots, N_h.
\end{equation}
To obtain a relation such as \eqref{eq:symbolic_SINDy}, we enforce it on the available data.
Denote by $\dot{\bs{X}}$ the matrix containing the corresponding time derivatives (computed numerically if unknown), and introduce $\bs{\Xi} = \{\bs{\xi}_1, \dots, \bs{\xi}_{N_h}\}\in\R^{r\times N_h}$, where the $k$-th column contains the nonzero coefficients for the equation $\dot{x}_k(t;\bs{\mu}) = f_k(t, \bs{x}(t;\bs{\mu});\bs{\mu})$.

We then compute $\bs{\Xi}$ by solving
\begin{equation}\label{eqSINDy}
    \dot{\bs{X}} = \bs{\Theta}(\bs{X}; \bs{\mu})\,\bs{\Xi},
\end{equation}
where $\bs{\Theta}(\bs{X};\bs{\mu})\in\R^{N_t 
N_\mu\times r}$ is the dictionary matrix containing the candidate function evaluations at the different time instances and parameter values.
Various algorithms exist to obtain a sparse solution $\bs{\Xi}$ to the regression problem in Equation~\eqref{eqSINDy}, for example, the \emph{LASSO} technique \cite{James_2013} and the \gls{stlsq} algorithm \cite{Brunton_2016}, with the latter being the one adopted in this work.

It is important to note that, in principle, this approach also accommodates time-dependent parameterized PDEs. Indeed, the PDE can be interpreted as a system of ordinary differential equations on the mesh points, and thus, the formulation in Equation \eqref{eqSINDyParam} remains valid, where $N_h$ now represents the number of degrees of freedom of the discretization. Essentially, the matrix $\bs{X}$ is built from the values of the numerical solution of the parameterized PDE, where each $\bs{x}(t_i,\bs{\mu}_j)\in\R^{N_h}$ corresponds to the state at time $t_i$ and parameter $\bs{\mu}_j$.

Applying \gls{sindy} to this system allows us to reproduce the time evolution at each degree of freedom, thereby reconstructing the full dynamics for each parameter value.

However, in many practical applications (especially when dealing with high-resolution discretizations or experimental data) $N_h$ may be extremely large—rendering direct solution computationally unfeasible. This difficulty is addressed in the next subsection.

\subsection{Reduction Techniques}
Directly applying \gls{sindy} to time-dependent PDEs may be prohibitively expensive. However, many high-dimensional systems exhibit underlying low-dimensional structures that can be discovered and exploited via dimensionality reduction techniques. Once a reduced set of coordinates is identified, \gls{sindy} can be applied on these new coordinates; the identified equations can be integrated and subsequently projected back to the original coordinates, ideally yielding an accurate description of the system.

To this end, many techniques have been employed in recent years, most notably the \gls{pod} (see, e.g., \cite{Quarteroni_2017}) and \textit{Dynamic Mode Decomposition} \cite{Jonathan_2024}, among others.
However, owing to the linearity of many such techniques, they may fall short when addressing nonlinear PDEs—especially time-dependent ones—since they might not capture essential nonlinear relationships in the data.
For these reasons, we adopt nonlinear reduction techniques based on autoencoder neural networks \cite{Goodfellow_2016}, as utilized in \cite{Champion_2019} and \cite{Conti_2023}. These architectures can identify a reduced set of coordinates via nonlinear transformations, potentially offering improved accuracy in capturing the system dynamics.

The remainder of this subsection details the tools used to develop an efficient \gls{rom} strategy: the \gls{pod} and the SINDy-AE architecture.

\subsubsection{The Proper Orthogonal Decomposition}\label{subsubsec:POD}
The Proper Orthogonal Decomposition (POD) is a powerful technique for deriving a low-dimensional representation of high-dimensional data. Intuitively, POD identifies the principal directions (modes) in the data that capture the maximum amount of variance, representing an optimal basis for the data in a linear reduced-dimensional space.

Consider a finite set $\{\bs{u}_j\}_{j=1}^{n_{\text{train}}}$ of snapshot vectors in $\R^{N_h}$, with $N_h$ the number of degrees of freedom in the discretization, that contains the solution information w.r.t.\ the time evolution and the parameter $\mu_j$. 
For any $N \leq n_{\text{train}}$, the POD basis of dimension $N$ extracted from $\{\bs{u}_j\}_{j=1}^{n_{\text{train}}}$ is defined as the set of $N$ orthonormal vectors solving the minimization problem:
\begin{equation}\label{eqPOD2}
\begin{cases}
\min\{E(\bs{z}_1, \ldots, \bs{z}_N) : \bs{z}_i \in \R^{N_h}, \, \bs{z}_i^T \bs{z}_j = \delta_{ij}, \, 1 \leq i,j \leq N\}\\[1mm]
\text{with} \quad E(\bs{z}_1, \ldots, \bs{z}_N) = \sum_{i=1}^{n_{\text{train}}} \| \bs{u}_i - \Pi_{Z_N} \bs{u}_i \|^2_2,
\end{cases}
\end{equation}
where $\Pi_{Z_N}$ denotes the projection onto the space spanned by $\{\bs{z}_1,\dots,\bs{z}_N\}$.
It can be shown that 
\begin{equation*}
E(\bs{\zeta}_1, \ldots, \bs{\zeta}_N) = \sum_{i=N+1}^{N_h} \sigma_i^2,
\end{equation*}
i.e., the approximation error using the POD basis equals the sum of the squares of the neglected singular values.
Thus, one may choose $\tilde{N}$ such that $E(\bs{\zeta}_1, \ldots, \bs{\zeta}_{\tilde{N}}) \leq \epsilon_{\text{tol}}^*$ for a prescribed tolerance. Equivalently, $\tilde{N}$ is chosen as the smallest $N$ such that
\begin{equation}\label{eq:POD_ENERGY}
I(N) = \frac{\sum_{i=1}^{N} \sigma_i^2}{\sum_{i=1}^{N_h} \sigma_i^2} \geq 1 - \delta,
\end{equation}
meaning that the relative energy contained in the first $\tilde{N}$ modes exceeds $1-\delta$ (with $\delta>0$ small). 

\subsubsection{Nested POD}\label{subsectNESTED_POD}
For time-dependent parameterized problems, performing a POD on the full snapshot matrix can be computationally prohibitive. To mitigate this, a nested (two-level) POD strategy is adopted \cite{nested_POD1, nested_POD2}.

First, for each parameter instance $ \mu_m \in \mathbb{P}_h = \{\mu_1, \mu_2, \dots, \mu_{N_\mu}\} \subset \mathbb{P} $, consider the snapshot matrix containing the evolution of the dynamical system
\[
\bs{S}^{N_h}_m = 
\begin{bmatrix}
\bs{u}^1_h(\mu_m) & \ldots & \bs{u}^{N_t}_h(\mu_m)
\end{bmatrix}
\in \mathbb{R}^{N_h \times N_t}.
\]
A first-level (local) POD is applied to $ \bs{S}^{N_h}_m $, yielding a set of $ K_m $ spatial modes. For each $\mu_m$, these modes are stored in the matrix $ \mathbb{V}_m \in \mathbb{R}^{N_h \times K_m}$, which are then concatenated to form 
\begin{equation}\label{eqNESTED_POD_global}
\mathbb{V}_{\text{global}} =
\begin{bmatrix}
\mathbb{V}_1 & \dots & \mathbb{V}_{N_\mu}
\end{bmatrix}
\in \mathbb{R}^{N_h \times K},
\end{equation}
where $ K = \sum_{m=1}^{N_\mu} K_m $. We remark that this procedure can be performed by setting a tolerance threshold on the energy retained by the spatial modes for each value of the parameter $\mu_m \in \mathbb{P}_h$, which can cause the number of bases $K_m$ extracted for each value of the parameter to be different.

Instead of directly using $ \mathbb{V}_{\text{global}} $  for projection, this matrix is treated as a snapshot matrix and a second-level (global) POD is performed. This additional POD yields an orthonormal basis $ \tilde{\mathbb{V}} \in \mathbb{R}^{N_h\times N_{\text{global}}} $ (with $ N_{\text{global}} \le K $), which is then employed in the projection of the full-order model.

\subsubsection{SINDy-AE-nested-POD architecture}
The SINDy-AE-nested-POD architecture combines the previously discussed methods to create an efficient reduced order model. 
To reduce computational complexity and obtain a more interpretable dynamical system, we perform nonlinear reduction through an \gls{ae}, identifying a new set of coordinates $\bs{z}(t)\in\R^n$ with $n\ll N_h$. 
The autoencoder network consists of a pair $(\bs{\varphi},\bs{\psi})$, where:
\begin{itemize}
    \item $\bs{\varphi}(\cdot) = \bs{\varphi}(\cdot; \bs{W}_\varphi): \R^{N_h}\to\R^n$ is the \textit{encoder}, parameterized by the weight matrix $\bs{W}_\varphi$,
    \item $\bs{\psi}(\cdot) = \bs{\psi}(\cdot; \bs{W}_\psi): \R^n\to\R^{N_h}$ is the \textit{decoder}, parameterized by the weight matrix $\bs{W}_\psi$,
\end{itemize}
where both $\bs{\varphi}$ and $\bs{\psi}$ are implemented as fully connected neural networks, and
the pair $(\bs{\varphi}, \bs{\psi})$ should satisfy:
\begin{equation}\label{eqAUTOENC} 
\bs{\psi}\left(\bs{\varphi}(\bs{y}; \bs{W}_\varphi); \bs{W}_\psi\right)\approx \bs{y}, \quad\forall\bs{y}\in\R^{N_h}.
\end{equation}
In these new coordinates $\bs{z}$ obtained through $\bs{\varphi}$, system \eqref{eqSINDyParam} becomes:
\begin{equation}\label{eqSINDyLATENT}
    \begin{cases}
   \dot{\bs{z}}(t; \bs{\mu}) = \bs{\tilde{f}}(t, \bs{z}(t; \bs{\mu}); \bs{\mu}), \\
   \bs{z}(t_0; \bs{\mu}) = \bs{z}_{0},
    \end{cases}
\end{equation}
where $\bs{\tilde{f}}(\cdot\,;\bs{\mu}):I\times\R^n\to\R^n$ describes the dynamics of the low-dimensional system with initial condition $(t_0, \bs{z}_0)\in I\times\R^n$, and $\bs{z}(t) = \bs{\varphi}(\bs{x}(t); \bs{W}_\varphi)$. 

The solution to Equation \eqref{eqSINDyLATENT} can then be computed and projected back to full coordinates using the decoder $\bs{\psi}$.
Clearly, the training of the network is a key part of the procedure and should be addressed with the right strategy to return the correct set of latent coordinates for the system identification. The procedure involves two distinct optimization problems. The first stems from the autoencoder network and the requirement that the encoder-decoder pair accurately approximates the identity mapping (see Equation \eqref{eqAUTOENC}). Denoting with $\bs{X}$ the snapshot matrix, this can be formulated as an optimization problem on the weights matrices $\bs{W}_\varphi$ and $\bs{W}_\psi$, respectively for the encoder $\bs{\varphi}$ and decoder $\bs{\psi}$ functions, as follows:
\begin{equation}\label{eqMINIMIZAUTO}
\argmin_{\bs{W}_\varphi, \bs{W}_\psi}\| \bs{\psi}\left(\bs{\varphi}(\bs{X};\,\bs{W}_\varphi);\,\bs{W}_\psi\right)-\bs{X}\|_2^2.
\end{equation}
It is important to note that $\bs{\varphi}$ (as well as $\bs{\psi}$) is applied row-wise, specifically:
\begin{equation*}
    \bs{\varphi}(\bs{X};\,\bs{W}_\varphi) = 
    \begin{bmatrix}
        \bs{\varphi}(\bs{x}(t_1; \bs{\mu}_1);\, \bs{W}_\varphi)^T \\
        \bs{\varphi}(\bs{x}(t_2; \bs{\mu}_1);\, \bs{W}_\varphi)^T\\
        \dots \\
        \bs{\varphi}(\bs{x}(t_{N_t}; \bs{\mu}_{N_\mu});\,
        \bs{W}_\varphi)^T
   \end{bmatrix}\in\R^{N_tN_\mu\times n}.
\end{equation*}
The second optimization problem relates to \gls{sindy} and is formulated in the $\bs{z}$ variables. Given a library of candidate functions $\bs{\Theta}$, this problem takes the form:
\begin{equation}\label{eqSINDyMinim2}
\argmin_{\bs{\Xi}} \|\dot{\bs{Z}}-\bs{\Theta}(\bs{Z};\bs{\mu})\bs{\Xi}  \|^2_2 + \lambda \|\bs{\Xi}\|_1,
\end{equation}
where $\lambda\in\R^+$ serves as a fixed regularization term to enforce sparsity of the solution, and $\bs{Z} = \bs{\varphi}(\bs{X}; \bs{W}_\varphi)$ represents the encoded snapshot matrix.

To capture the most important dynamical features while retaining interpretable terms, the optimization tasks \eqref{eqMINIMIZAUTO} and~\eqref{eqSINDyMinim2} must be addressed simultaneously through a unified neural network architecture that incorporates both tasks into a single optimization problem. The sparse regression terms are inserted as additional regularization components in the \gls{ae} loss function, with the sparse coefficients $\bs{\Xi}$ treated as network variables and estimated concurrently with $\bs{\varphi}$ and $\bs{\psi}$ during training.

The resulting optimization problem is given by:
\begin{equation}\label{eqSINDyAEMinim}
\begin{split}
\underset{\bs{W}_\varphi, \bs{W}_\psi, \bs{\Xi}}{\argmin} \Bigg(\underbrace{\|\bs{X} - \bs{\psi}(\bs{\varphi}(\bs{X};\bs{W}_\varphi);\bs{W}_\psi)\|^2_{2}}_{\text{Autoencoder loss}}+ \underbrace{\lambda_1 \|\dot{\bs{Z}} - \bs{\Theta}(\bs{Z};\bs{\mu})\bs{\Xi}\|^2_{2}+\lambda_2\|\bs{\Xi}\|_{1}}_{\text{Sparse regression loss}} \\+ \underbrace{\lambda_3 \|\dot{\bs{X}} - \nabla_{\bs{z}}\psi(\bs{Z};\bs{W}_\psi) \bs{\Theta}(\bs{Z};\bs{\mu})\bs{\Xi}\|_{2}^{2}}_{\text{Consistency loss}} \Bigg)\,.
\end{split}
\end{equation}
Recalling that $\dot{\bs{z}}(t) = \dot{\bs{\varphi}}(\bs{x}(t)) = \nabla_{\bs{x}}\bs{\varphi}(\bs{x}(t))\dot{\bs{x}}(t)$ by the chain rule, the consistency loss term in \eqref{eqSINDyAEMinim} ensures that the time derivative of the network output obtained via \gls{sindy} matches the input time derivatives $\dot{\bs{X}}$. The coefficients $\lambda_1,\, \lambda_2,\, \lambda_3\in\R^+$ are hyperparameters of the technique.
These values must be adjusted according to the specific problem and dataset.

The latent dimension $n$ is another crucial hyperparameter requiring careful tuning. When the dimension of the low-dimensional manifold is known a priori, $n$ can be set accordingly. Otherwise, it can be determined by analyzing the behavior of the SINDy-AE architecture's loss function with respect to $n$, selecting the value beyond which further dimension reduction leads to significant error increases \cite{Conti_2023}. While SINDy-AE can operate with $n$ larger than the manifold dimension, this may introduce computational challenges, as \gls{sindy}'s efficiency rapidly decreases with increasing $n$, potentially complicating the dynamics identification task.

The SINDy-AE-nested-POD architecture operates in three primary stages. First, nested POD efficiently computes a reduced basis for the high-dimensional data, handling temporal and parametric variations separately. Second, the autoencoder further reduces dimensionality through nonlinear transformations, capturing complex relationships that linear methods might miss. Finally, SINDy identifies sparse, interpretable governing equations in the latent space for prediction and analysis.

\begin{figure}[ht]
    \centering
    \includegraphics[width=0.65\textwidth]{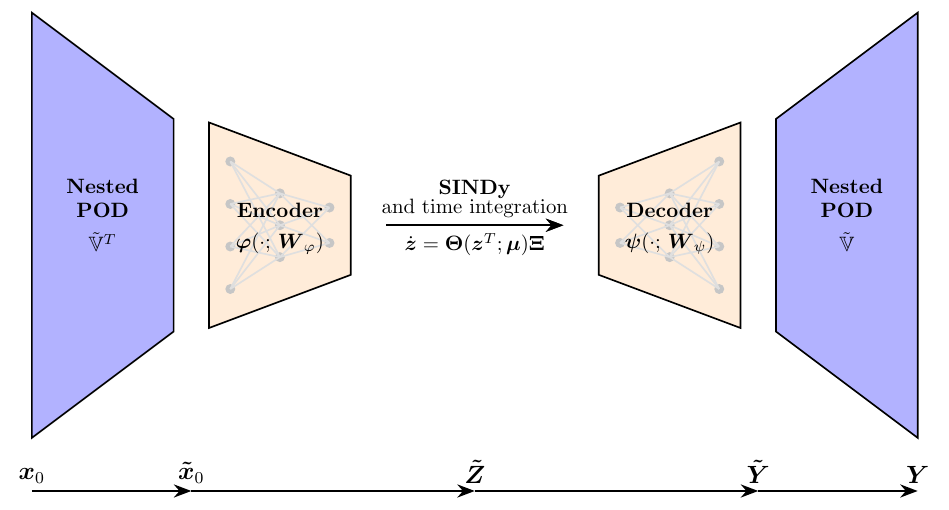}
    \caption{Representation of the SINDy-AE-nested-POD architecture.}
    \label{figPOD_AE_SINDy}
\end{figure}

The \textit{offline training} phase consists of minimizing the loss function through backpropagation using the ADAM algorithm. Once trained, the model can be queried \textit{online} to reconstruct the complete time evolution of the full-order system for different initial conditions and new parameter instances.
The procedure, depicted in Figure \ref{figPOD_AE_SINDy} and
summarized in Algorithm \ref{alg:online_sindy_ae}, reads as follows:
 \begin{enumerate}[label=(\roman*)]
    \item\label{ItemInitCondCompress} \textbf{Encode} the initial condition $\bs{x}_0$ of the dynamical system \eqref{eqSINDyParam} by first mapping it to \gls{pod} coordinates and then computing its latent representation, obtaining $\bs{z}_0 = \bs{\varphi}(\tilde{\mathbb{V}}^T\bs{x}_0;\, \bs{W}_\varphi)$;
    
     \item\label{ItemLatentIntegration} \textbf{Integrate} Equation~\eqref{eqSINDyLATENT}, obtained exploiting \gls{sindy} on latent variables, using a time-marching technique, such as Runge-Kutta or multistep methods.
     That is, given $\bs{z}_0$, recover the latent dynamics by computing $\bs{z}_1, \dots, \bs{z}_{N_\text{end}}$ as approximations of the solution $\bs{z}(t)$ of~\eqref{eqSINDyLATENT} at the $N_\text{end}+1$ points used for time discretization $[t_0,\, T_\text{end}]$;
     
     \item\label{ItemDecompression} \textbf{Decode} the latent solution using the decoder, project back onto physical coordinates, and reconstruct the full-order solution via $\tilde{\mathbb{V}}\bs{\psi}(\bs{z}_i\,; \bs{W}_\psi)\approx \bs{x}(t_i),$ for $i = 0, \dots, N_{\text{end}}$.
 \end{enumerate}

This online procedure is computationally efficient since the encoding and decoding steps merely involve evaluating pre-trained \gls{nn} functions, while the integration step operates on a low-dimensional system. Due to \gls{sindy}'s ability to extrapolate in time and parameters, $T_{\text{end}}$ can significantly exceed the final time used for model training, performing an extrapolation task. Furthermore, when only the steady state of \eqref{eqParamDyn} is of interest, one can integrate the latent system \eqref{eqSINDyLATENT} to its steady solution and then decode using $\bs{\psi}$, efficiently obtaining the desired result.

\begin{algorithm}
\SetAlgoNlRelativeSize{0}
\caption{Pseudo-code for the online phase}
\label{alg:online_sindy_ae}
\KwIn{$\bs{x}_0$ initial condition, $\mu_\text{test}$ parameter value}
\KwOut{$(\bs{x}_0, \dots, \bs{x}_{N_\text{end}})$ simulated trajectory}
$\tilde{\bs{x}}_0 = \tilde{\mathbb{V}}^T\bs{x}_0$\tcp*{POD projection}
$\bs{z}_0 = \bs{\varphi}(\tilde{\bs{x}}_0^T;\, \bs{W}_\varphi)$\tcp*{Nonlinear encoding}
$(\bs{z}_1, \dots, \bs{z}_{N_\text{end}})= \text{TimeStepping}(\tilde{\bs{f}}(\cdot;\, \mu_\text{test}), \bs{z}_0)$\tcp*{Integration in time}

\For{$i=0, \dots, N_{\text{end}}$
}{
$\tilde{\bs{x}}_i = \bs{\psi}(\bs{z}_i\,; \bs{W}_\psi)$\tcp*{Nonlinear decoding}
$\bs{x}_i = \tilde{\mathbb{V}}\tilde{\bs{x}}_i$\tcp*{POD reconstruction}
}
\end{algorithm}

\begin{remark}[SINDy-AE with nested POD]\label{rmkNESTED}
As discussed previously, a further enhancement of the plain SINDy-AE architecture toward an efficient \gls{rom} has been proposed in \cite{Conti_2023}, introducing the \gls{pod} projection as a reduction technique. 
By projecting the data using the matrix $\tilde{\mathbb{V}} = [\bs{\zeta}_1 \dots \bs{\zeta}_{N_\text{POD}}]\in\R^{N_h\times N_\text{POD}}$ formed by the first $N_\text{POD}$ dominant modes,
this strategy enables a more efficient encoding process, as the autoencoder now only needs to reduce the dimension from a relatively low number $N_\text{POD}$ of features to $n$, rather than directly handling the full $N_h$ degrees of freedom. Reducing a reduced set of variables could seem counterintuitive, but it could prove to be fundamental when the system exhibits slow Kolmogorov $n$-width decay, e.g.\ advection-dominated problems, for which a large number of modes are necessary for the linear compression.

In this work, we employ nested \gls{pod} as the preliminary reduction technique instead of standard POD. The fundamental procedure remains unchanged; the only modification is the use of matrix $\tilde{\mathbb{V}}$ resulting from the nested POD for projections, enabling a more efficient reduction when working with large datasets containing many degrees of freedom.
\end{remark}

\section{Numerical results}\label{sec:numerical}
In this section, we present numerical results for the two test cases introduced in Sections \ref{subsec:coanda_problem} and \ref{subsec:hopf_problem}. Since problem settings, boundary conditions, and numerical discretization schemes have been thoroughly described in their respective sections, we focus here on demonstrating the effectiveness of the SINDy-AE-nested-POD architecture in capturing the bifurcating phenomena. For each case, we analyze the architecture's ability to learn the underlying dynamics, reconstruct the bifurcation diagrams, and extrapolate solutions in both parameter and time domains. We also investigate the sparsity of the identified governing equations through a dedicated study of the nonparametric case. The results showcase the methodology's capability to accurately reproduce complex dynamical behaviors while maintaining computational efficiency.
The network architecture and training parameters for the proposed SINDy-AE models are summarized in Table \ref{tab:network_params}.

The code for data collection in the pitchfork bifurcation case was developed using the FEniCS library \cite{LoggEtal_2012}, whereas for the Hopf bifurcation we exploited the OpenFOAM library \cite{Weller_1998}. Furthermore, the \gls{sindy}-\gls{ae}-nested-POD architecture was developed on top of the publicly available \gls{sindy}-\gls{ae} implementation provided in the GitHub repository \cite{VVV_code}.
In particular, to perform the nested POD we used the OpenFOAM-based ITHACA-FV library \cite{Stabile_2018,Stabile2017CAIM}.

\begin{table}[htbp]
    \centering
    \begin{tabular}{llll}
        \hline
        Parameter & Symmetry-breaking & Hopf bifurcation & Hopf (nonparametric) \\
        \hline
        Encoder layers & [32, 8, 4] & [64, 32, 16, 8, 4, 3] & [64, 32, 16, 8, 4] \\
        Decoder layers & [4, 8, 32] & [3, 4, 8, 16, 32, 64] & [4, 8, 16, 32, 64] \\
        Latent dimension & 2 & 2 & 2 \\
        Training epochs & 5000 & 350 & 10000 \\
        Learning rate & $10^{-5}$ & $10^{-5}$ & $10^{-5}$ \\
        Batch size & 64 & 64 & 64 \\
        $\lambda_1$ & $10^{-10}$ & $10^{-9}$ & $10^{-9}$ \\
        $\lambda_2$ & $10^{-6}$ & 0 & 0 \\
        $\lambda_3$ & 0 & $10^{-10}$ & $3\cdot 10^{-9}$ \\
        \hline
    \end{tabular}
    \caption{Network architecture and training parameters for both test cases.}
    \label{tab:network_params}
\end{table}

\subsection{Symmetry-breaking bifurcation}\label{subsec:coanda_results}
Given the previous analysis, we considered a discretization of $\mathbb{P}=[0.75, 1.05]$ comprising $20$ equispaced values of $\mu$.
We performed the nested \gls{pod} as the preliminary dimensionality reduction technique for the \gls{sindy}-AE architecture.
Specifically, we computed local \glspl{pod} on the snapshot matrices corresponding to each parameter value, with $K_m$ modes chosen for each $\mu_m$ according to the criterion in Equation \eqref{eq:POD_ENERGY}, using $\delta_\text{local}=10^{-6}$.
Similarly, $N_\text{global}$ was determined using the same criterion using $\delta_\text{global}=10^{-5}$.
The projection onto \gls{pod} coordinates was then performed for the full-order snapshots corresponding to the selected values of $\mu$.
After this operation, the \gls{pod} coefficients of $u_1$, $u_2$, and $p$ were independently scaled to have zero mean and unit variance to facilitate network training, then stacked and passed as input to the architecture.

Since we are not interested in reproducing the initial transient phase, as it does not convey useful information about the bifurcation diagram, we considered only snapshots corresponding to time instances in the interval $[8.0, 70.0]$.
For parameter values where the simulation stopped before reaching $T_\text{end}=70.0$, we duplicated the last snapshot until the final time $T_\text{end}$ (they are all approximations of the steady-state solution).
We used $90\%$ of the computed snapshots for training the network and \gls{sindy}, reserving the remaining $10\%$ to test the model's ability to extrapolate in time.
To ensure better stability, we trained a separate \gls{sindy} model on the latent variables using the Python library \emph{pysindy} \cite{PYSINDy1, PYSINDy2}.
The library was constructed using a polynomial basis of degree $2$ with respect to both state variables and parameter, employing both ensembling and library ensembling to improve robustness to noise \cite{SINDy_ensemble}.
The STLSQ threshold $\tau$ was set to $0.01$.

\begin{figure}[htbp]
\centering
\begin{subfigure}{\linewidth}
\centering
\includegraphics[width=0.9\linewidth]{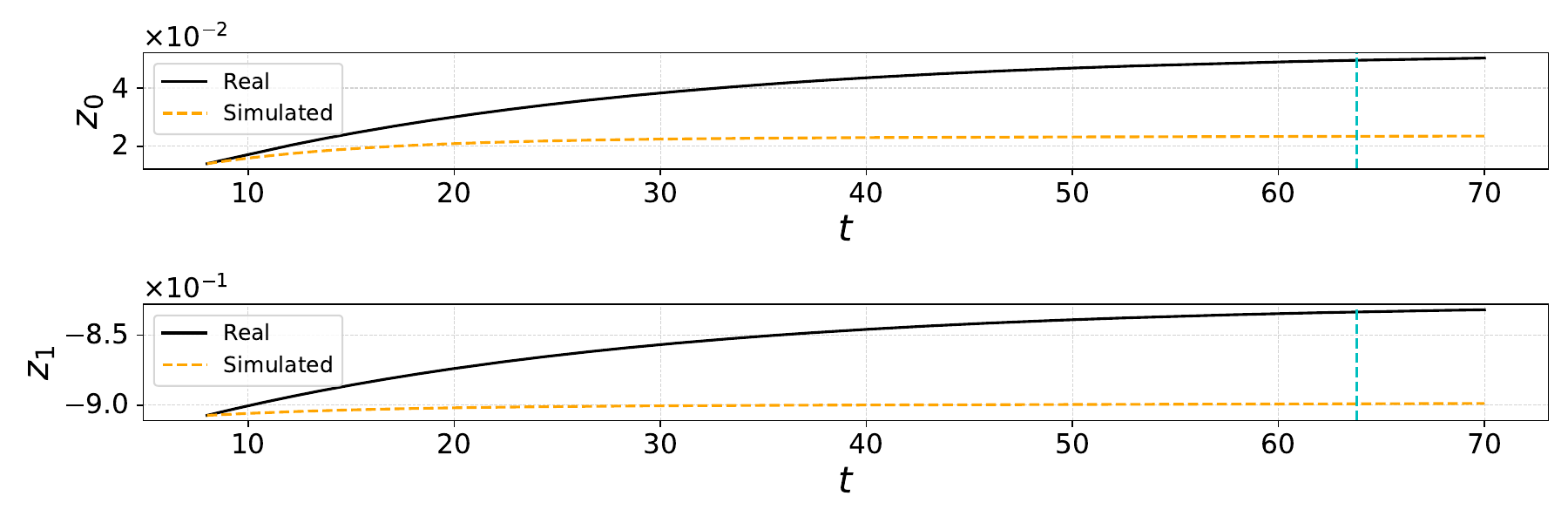}
\caption{Latent trajectory for $\mu=1.01$}
\label{latentMU_PREBIF}
\end{subfigure}

\vspace{0.5cm}

\begin{subfigure}{\linewidth}
\centering
\includegraphics[width=0.9\linewidth]{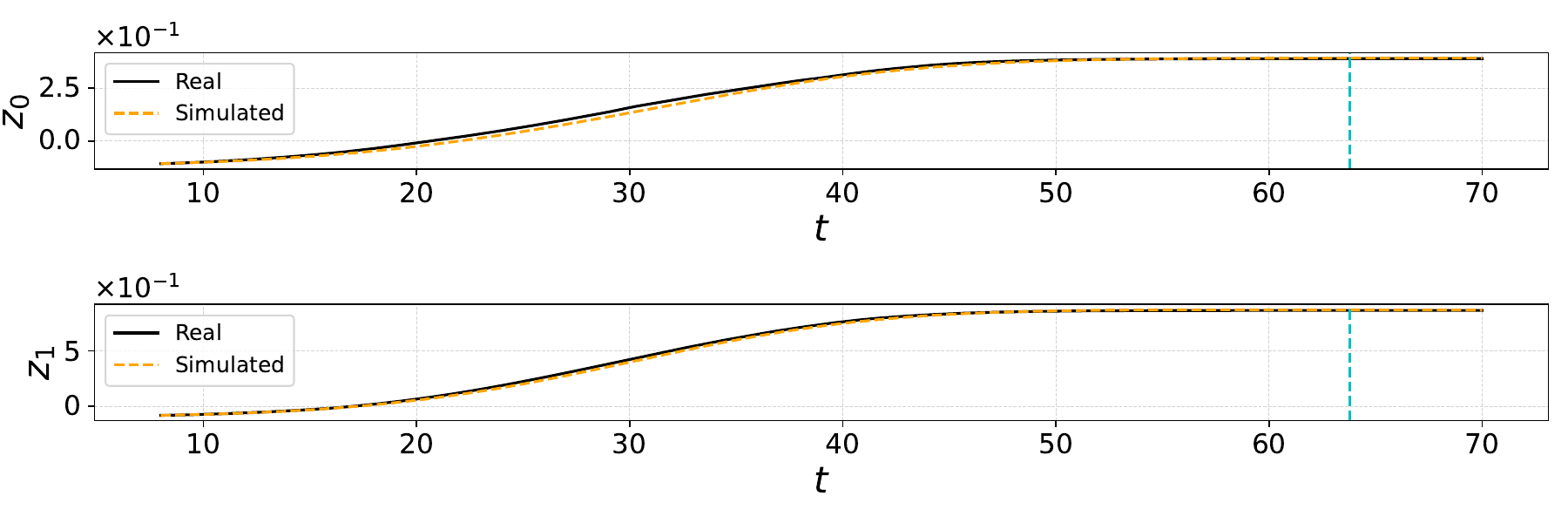}
\caption{Latent trajectory for $\mu=0.81$}
\label{latentMU_BIF}
\end{subfigure}
\caption{Evolution of the latent trajectories for two testing values of $\mu$.}
\label{fig:latent_coanda}
\end{figure}

Figure \ref{fig:latent_coanda} shows the real and simulated latent trajectories for the test values $\mu = 0.81$ and $\mu = 1.01$.
We observe a general agreement of the real and simulated evolutions. SINDy reconstruction in the non-bifurcating regime \ref{latentMU_PREBIF} may appear less accurate, but the error magnitude is lower than the evolution scale in \ref{latentMU_BIF}, indicating that the model correctly identifies the ``flat" dynamics converging to the symmetric steady state. The dashed blue line denotes the end of the training interval.

\begin{figure}[htbp]
\begin{subfigure}{0.49\linewidth}
\centering
\includegraphics[width=\linewidth]{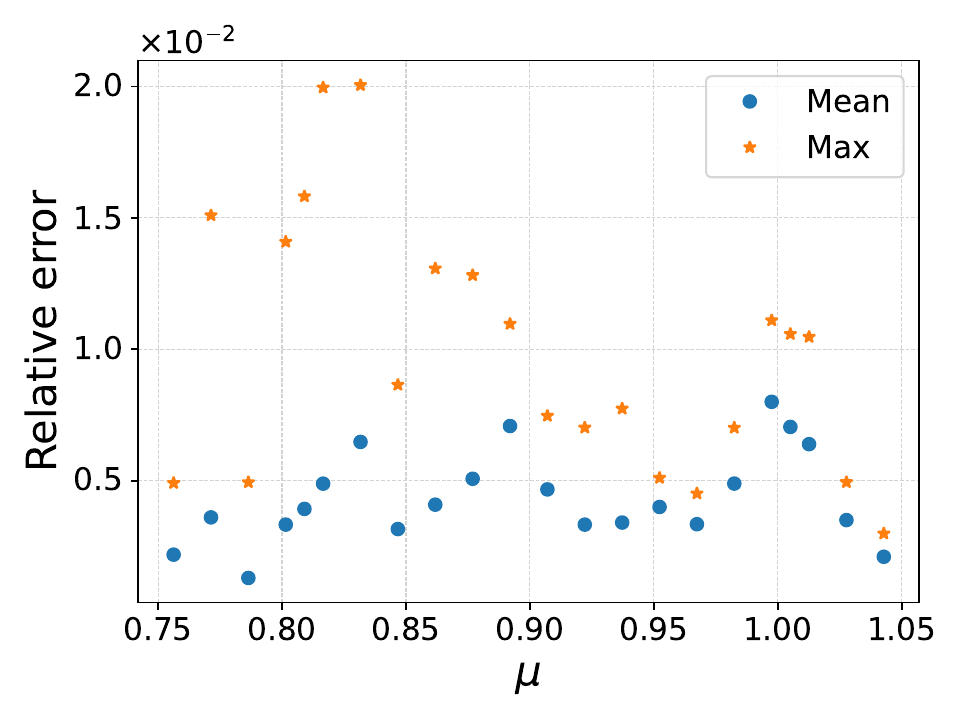}
\caption{Errors on $u_1$}
\label{fig:err_u1}
\end{subfigure}
\hfill
\begin{subfigure}{0.49\linewidth}
\centering
\includegraphics[width=\linewidth]{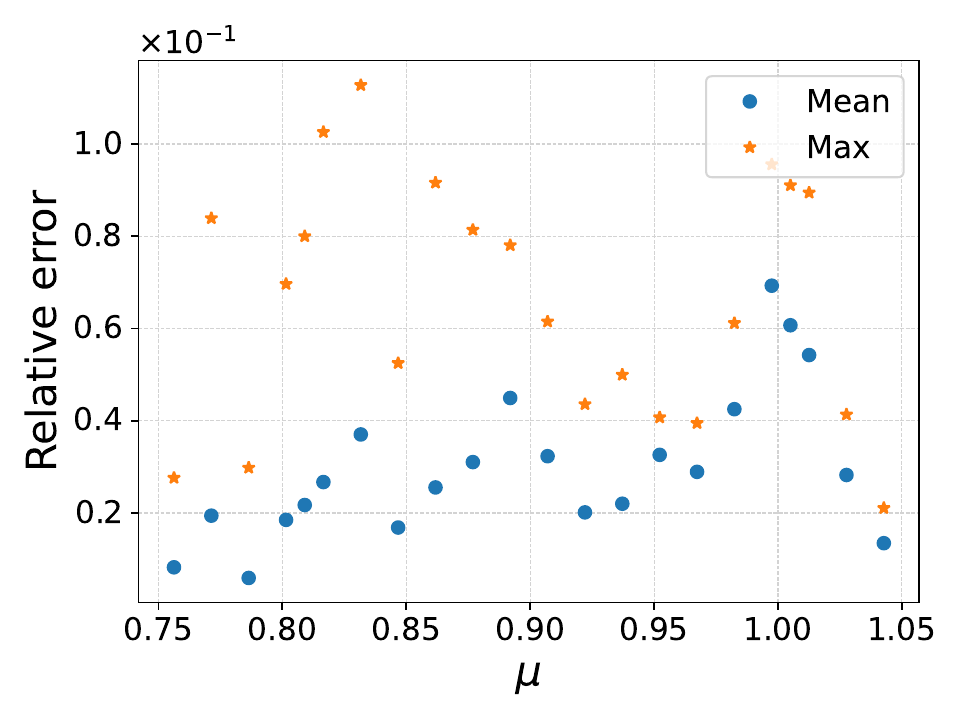}
\caption{Errors on $u_2$}
\end{subfigure}
\caption{Mean and maximum errors of the SINDy-AE-nested-POD architecture against the values of $\mu\in\mathbb{P}_h$.}
\label{fig:coanda_error}
\end{figure}

In Figure \ref{fig:coanda_error} the mean and maximum errors introduced by the SINDy-AE-nested-POD architecture are plotted against the values of $\mu\in\mathbb{P}_h$. We have not reported the errors for $p$ due to their similarity to those for $u_1$ in Figure \ref{fig:err_u1}. As it is clear from the plots, the error for the vertical components of the velocity $u_2$ is one order of magnitude higher than its horizontal counterpart. Indeed, the symmetry-breaking phenomenon is essentially related to the vertical axis, and being $u_2$ the primary driver of the bifurcating behavior, such complexity causes the approximation of the latent behavior to be less robust.

\begin{figure}[htbp]
\centering 
\begin{subfigure}{0.47\linewidth}
\centering
\includegraphics[width=\linewidth, trim=0 10 0 130, clip]{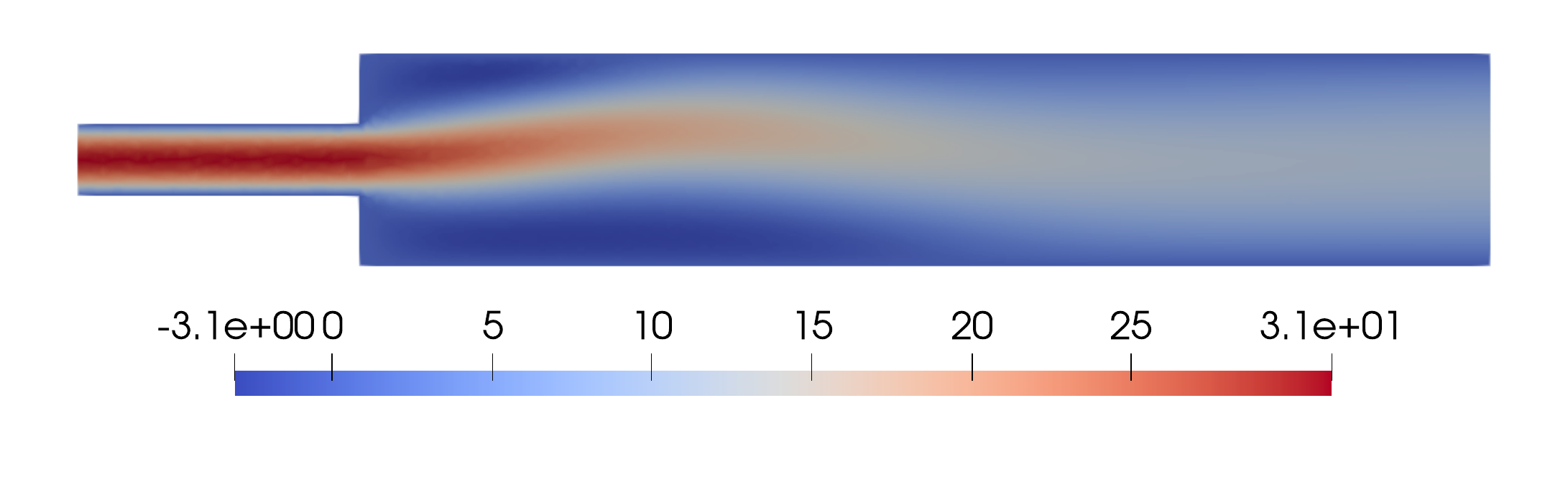}
\caption{Exact}
\end{subfigure}
\begin{subfigure}{0.47\linewidth}
\centering
\includegraphics[width=\linewidth, trim=0 10 0 130, clip]{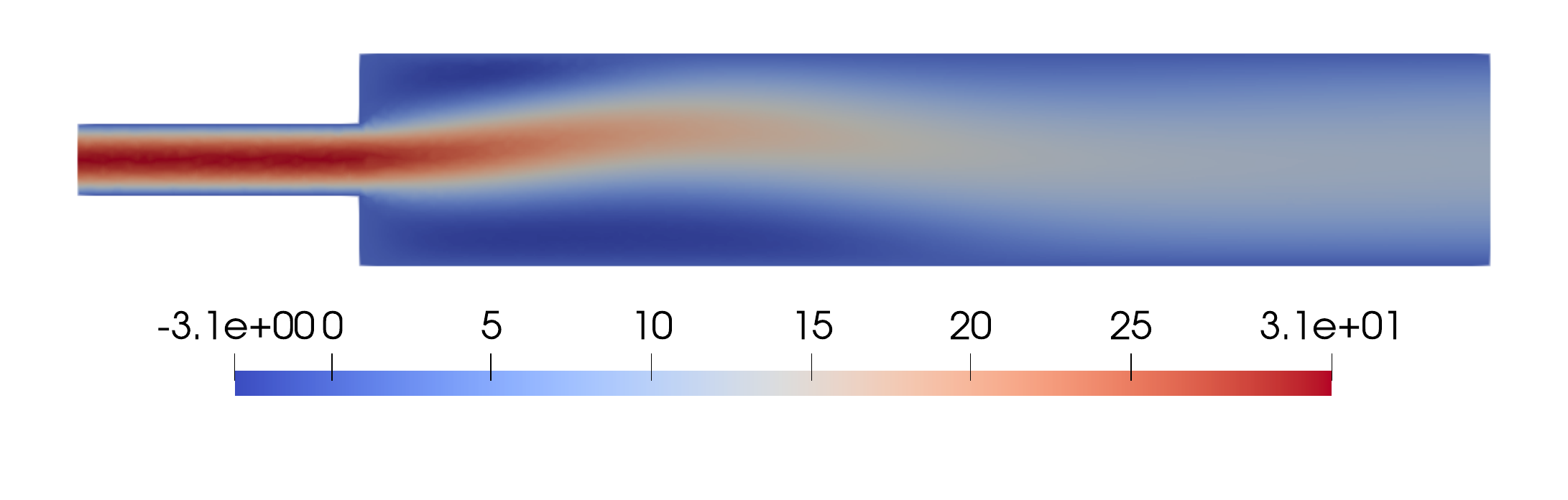}
\caption{Simulated}
\end{subfigure}
\begin{subfigure}{0.47\linewidth}
\centering
\includegraphics[width=\linewidth, trim=0 10 0 130, clip]{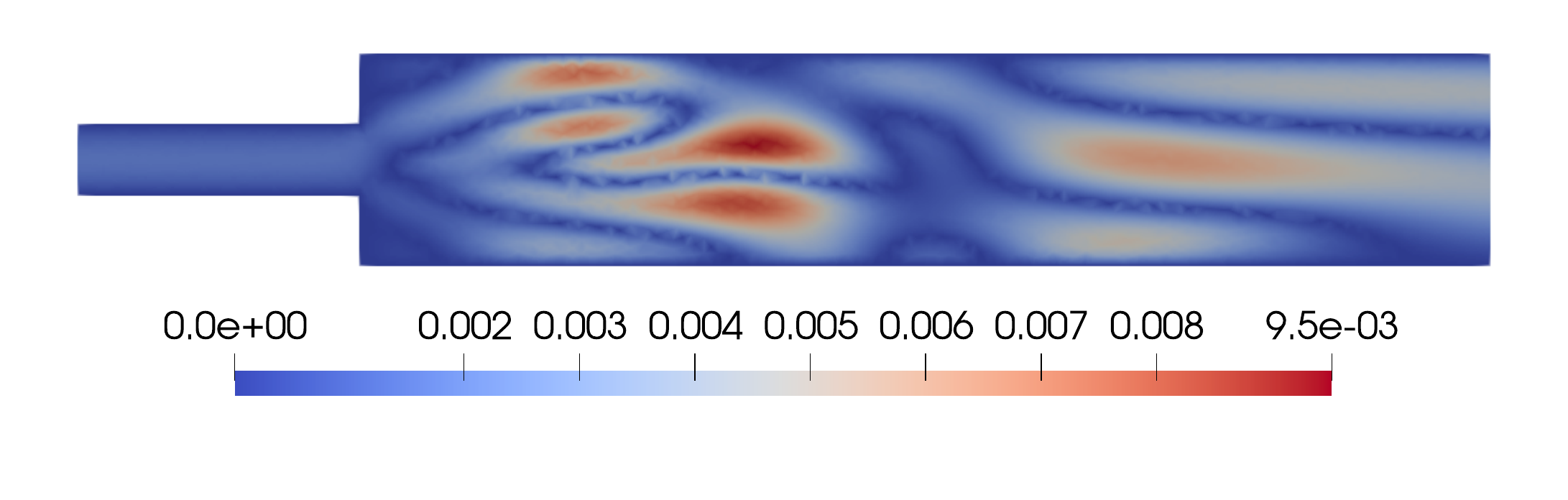}
\caption{Error}
\end{subfigure}
\caption{Comparison between the exact and simulated horizontal component of the velocity field $u_1$, and the corresponding error field for $\mu=0.81$ at $t=69.81$.}
\label{fig:coanda_reconstruction}
\end{figure}

We show in Figure \ref{fig:coanda_reconstruction} a comparison between a high-fidelity snapshot and its reconstruction, along with the resulting error field. We observe that despite the complexities in approximating the wall-hugging latent dynamics, the methodology is capable of correctly reconstructing the bifurcating steady states for an unseen value of the parameter. 

\begin{figure}[htbp]
    \centering
    \includegraphics[width=0.5\textwidth]{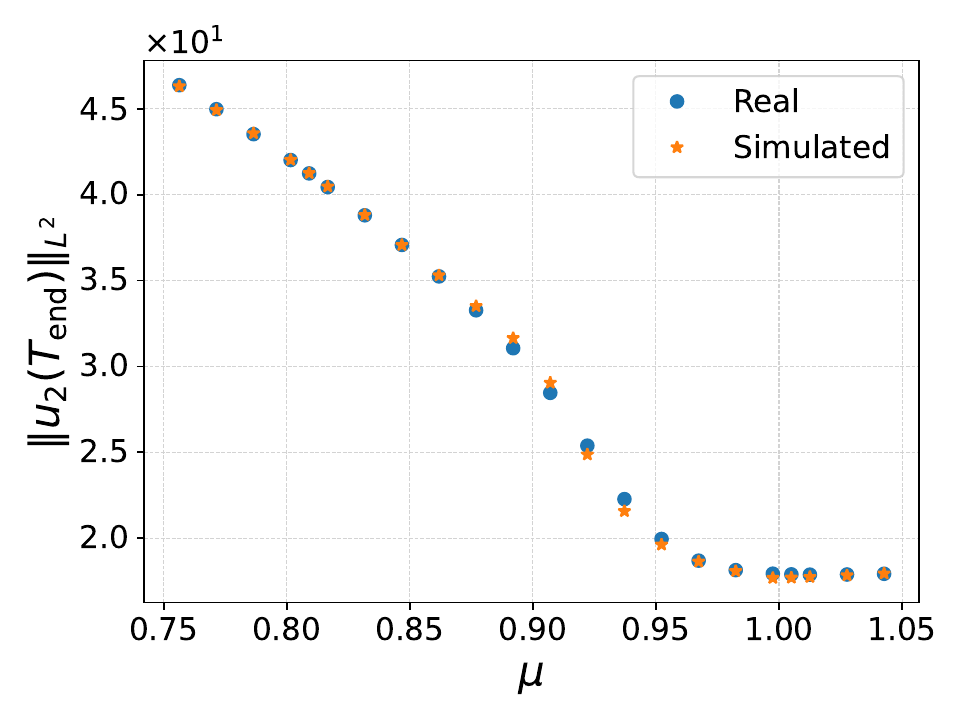}
    \caption{Real and simulated bifurcation diagram in terms of the $L^2$-norm of the vertical velocity $u_2$.}
    \label{fig:coanda_diagram}
\end{figure}

Finally, Figure \ref{fig:coanda_diagram} shows the real and simulated bifurcation diagram exploiting the $L^2$-norm of the vertical velocity $u_2$ as the scalar quantity of interest.
Once again, we note that despite the complex evolution of $u_2$, the methodology successfully reconstructs the bifurcation diagram with good accuracy.

\subsection{Hopf bifurcation}
For this case, we also employed the nested \gls{pod} as the preliminary dimensionality reduction technique for the \gls{sindy}-AE architecture.
Specifically, we performed POD on the snapshot matrices corresponding to each parameter value, setting $K_m=100$, while choosing $N_\text{global}$ as $32$ for both velocity and pressure.

Following the same preprocessing strategy, the POD coefficients of $\bs{u}$ and $p$ were independently scaled to have zero mean and unit variance to facilitate network training, then stacked and passed as input to the architecture\footnote{In the previous case, the vertical component of the velocity $u_2$ was much harder to approximate than $u_1$, hence suggesting to treat them separately. Here, instead, we stacked the matrices containing the snapshots of $u_1$ and $u_2$, performing the POD directly on the snapshot matrix $\bs{u}\in\R^{2N_h}$.}.

Among the $40$ equispaced values of $Re\in\mathbb{P}=[346, 356]$, for which we collected the full-order solution, we used $38$ for training and reserved two ($Re=347.28$ and $Re=354.97$) for testing.

To ensure a fully developed flow regime and focus on the steady-unsteady transition, we considered the system's evolution over the time interval $I=[1900, 2500]$, avoiding the initial diffusive transient.
Since the time step between successive snapshots was $1$, we interpolated the POD coefficients matrix using cubic natural spline functions and evaluated them at intervals of $dt=0.01$ seconds to ensure smooth computation of derivatives for SINDy.

In this case as well, we trained a separate \gls{sindy} model on the latent variables using pysindy.
The model was trained on a subset of the initial time interval $I_\text{train}=[2275, 2477]$ for efficiency purposes (using $dt=0.5\times10^{-2}$).
The library was constructed using polynomials of degree $3$ with respect to state variables and degree $2$ with respect to the parameter.
The \gls{stlsq} threshold $\tau$ in was set to $0.01$.

\begin{remark}[Error analysis]

Figure \ref{figERR_NESTED_POD} shows the mean projection error onto the nested POD basis for all values of $Re\in \mathbb{P}_h$. The projection error remains relatively low (on the order of $10^{-2}$) across the Reynolds number range, but shows a notable increase around $Re\approx 350.6$, which corresponds to the onset of the Hopf bifurcation. Before this critical point, the error for both the velocity and pressure fields remain constant. After the bifurcation point, both fields exhibit increasing errors, with pressure showing a more pronounced rise. This behavior is expected, as the flow dynamics becomes more complex after the bifurcation, making the projection onto a fixed number of POD modes slightly less accurate.

\begin{figure}[htbp]
\centering
\includegraphics[width=0.48\linewidth]{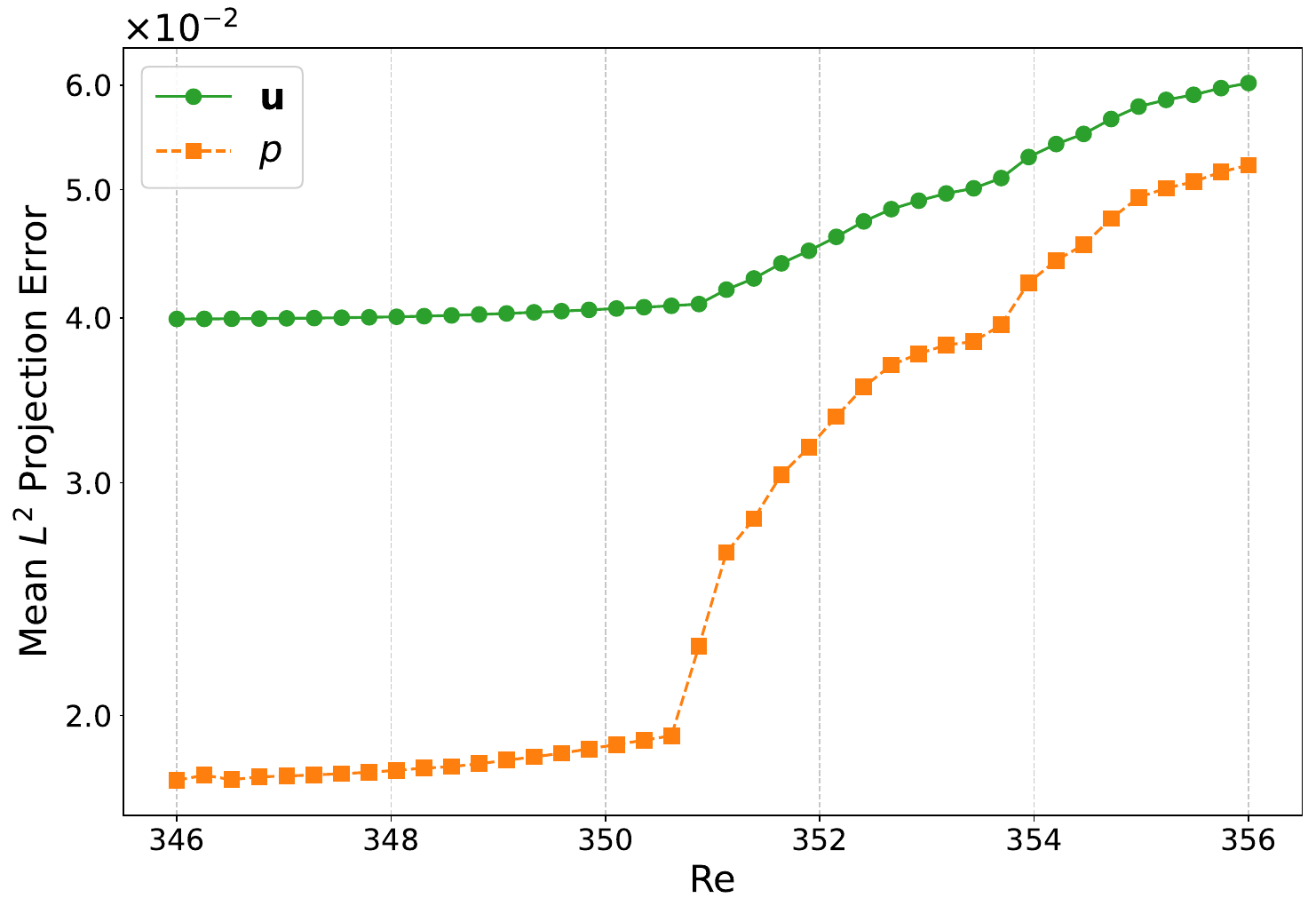}
\caption{Mean relative projection error for $\bs{u}$ and $p$ using nested POD for $Re\in\mathbb{P}_h$ with $N=32$ basis.}
\label{figERR_NESTED_POD}
\end{figure}

The relatively low error introduced by the \gls{pod} projection with $N_\text{POD}=32$ indicates that a meaningful analysis of the system can be conducted using the \gls{pod} coefficients alone, as they provide a sufficiently accurate representation of the full-order trajectories. Therefore, in the following numerical results, we focus exclusively on the \gls{pod} coefficients, as they are sufficient to recover the bifurcation diagram.
\end{remark}

\begin{figure}[htbp]
\centering
\begin{subfigure}{0.9\linewidth}
\centering
\includegraphics[width=\linewidth]{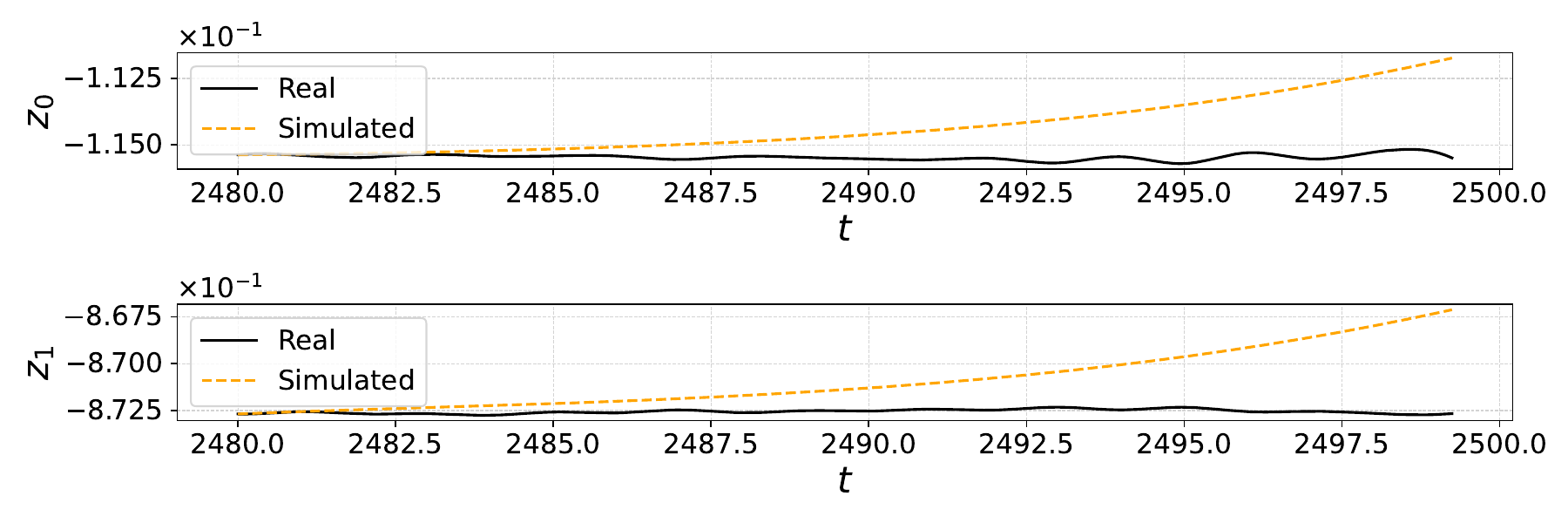}
\caption{Latent trajectory for $Re=347.28$}
\label{latentRE_PREBIF}
\end{subfigure}

\vspace{0.5cm}

\begin{subfigure}{0.9\linewidth}
\centering
\includegraphics[width=\linewidth]{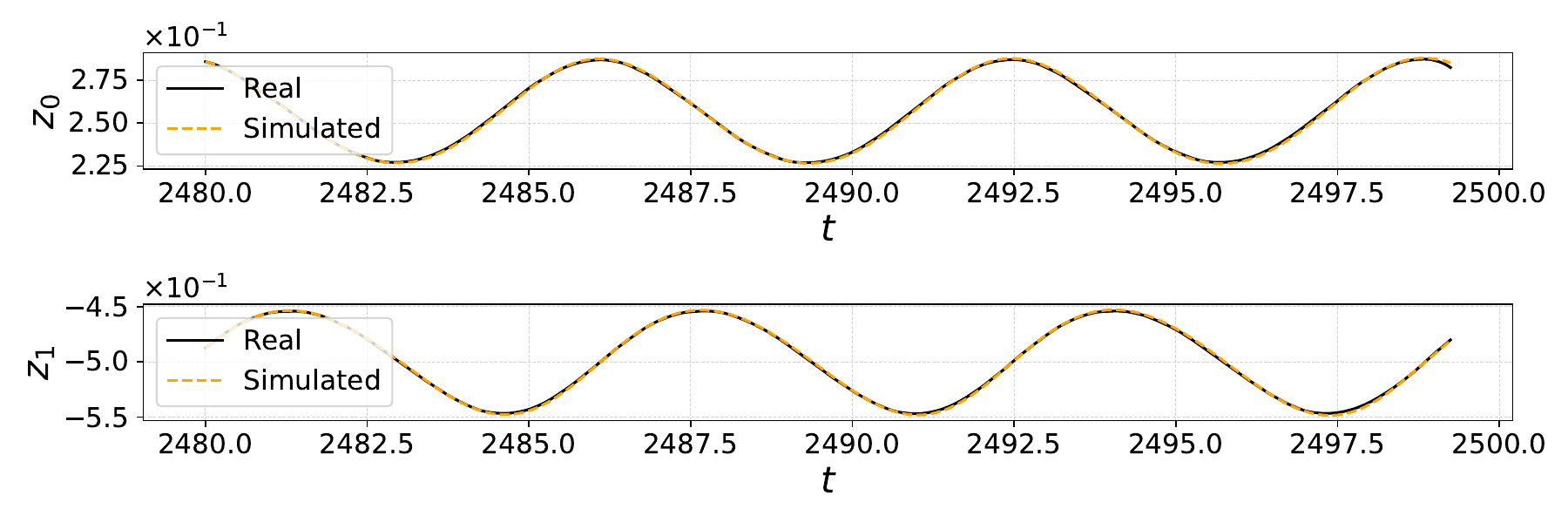}
\caption{Latent trajectory for $Re=354.97$}
\label{latentRE_BIF}
\end{subfigure}
\caption{Evolution of the latent trajectories for two testing values of $Re$.}
\label{figQUAINI_LATENT_P}
\end{figure}

To test the model's ability to extrapolate in time, we integrated over $[2480, 2500]$ (not seen during training) for each value of $Re$, using the POD coefficients at $t=2480$ as initial conditions.
Figure \ref{figQUAINI_LATENT_P} shows the real and simulated latent trajectories for the two test values of $Re$.

Similarly to the previous test case, the simulated evolution in \ref{latentRE_PREBIF} appears to deviate from the original one, but the error magnitude is much smaller than the oscillation amplitude in Figure \ref{latentRE_BIF}, indicating that SINDy correctly identifies the pre-bifurcation dynamics as nearly constant. On the other hand, we obtain very accurate results in the extrapolatory regime in parameter and time spaces, correctly identifying the periodic behavior of the Hopf bifurcation also at the reduced level. 

\begin{figure}[htbp]
\begin{subfigure}{0.49\linewidth}
\centering
\includegraphics[width=\linewidth]{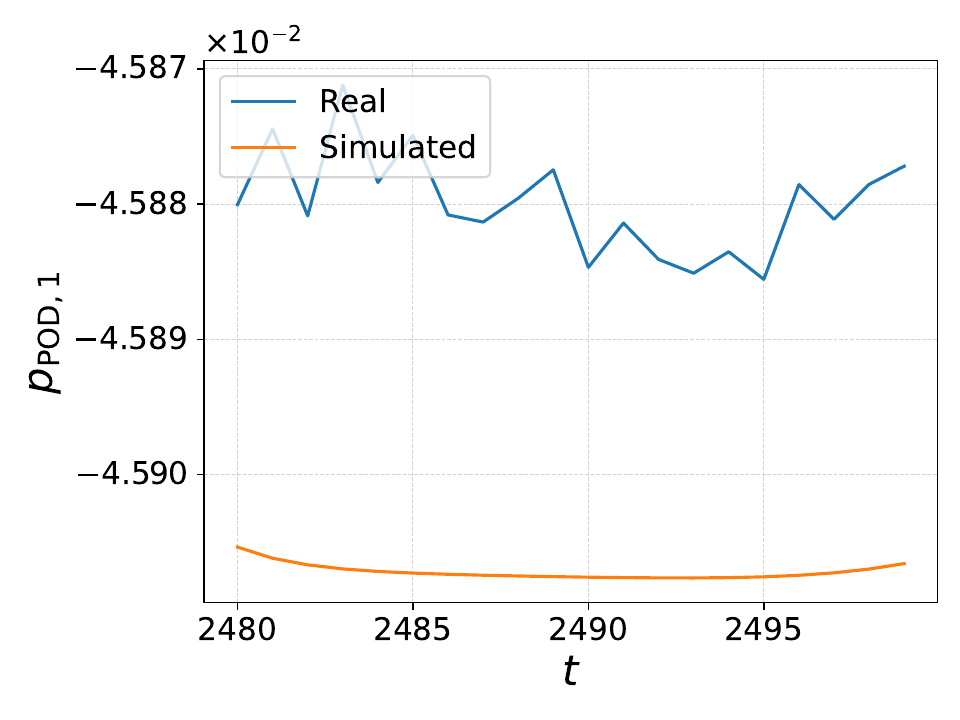}
\caption{$Re=347.28$}
\label{figPARAM_P_RECONSTR_PREBIF}
\end{subfigure}
\hfill
\begin{subfigure}{0.49\linewidth}
\centering
\includegraphics[width=\linewidth]{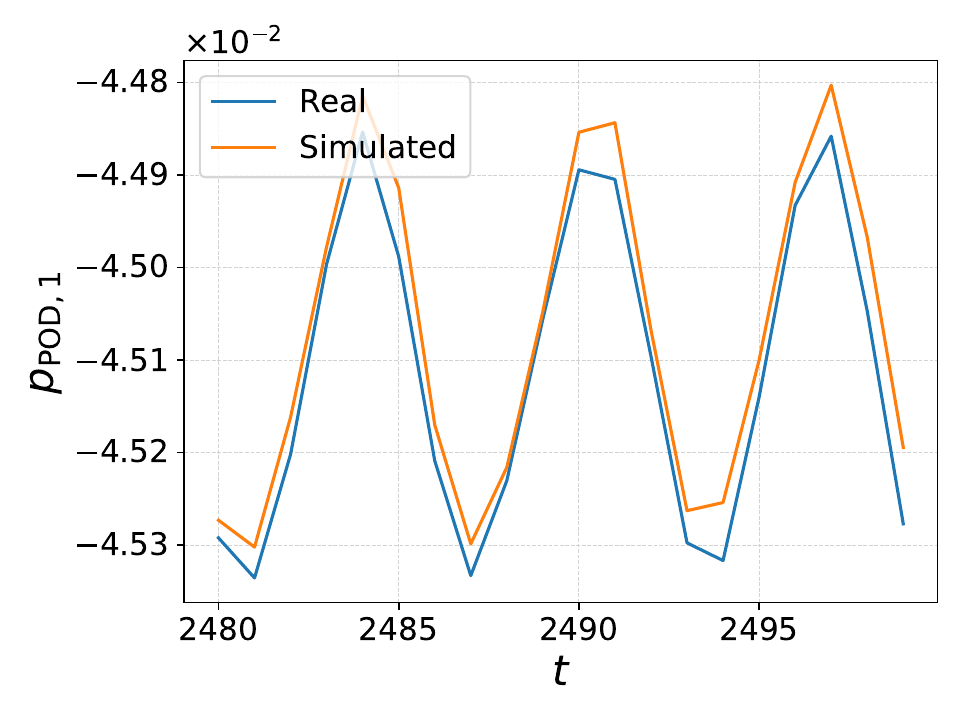}
\caption{$Re=354.97$}
\label{figPARAM_P_RECONSTR_BIF}
\end{subfigure}
\caption{Evolution of the first POD coefficient of $p$ for two testing values of $Re$.}
\label{figQUAINI_P_TRAJ}
\end{figure}

Figure \ref{figQUAINI_P_TRAJ} shows the real and simulated evolution of the first POD coefficients of $p$, from which we can observe the complete change in the dynamics of the system transitioning from steady state to unsteady periodic behavior.
We omit the evolution of the first POD coefficients of $\bs{u}$ as they exhibit behavior similar to that of $p$.

\begin{figure}[htbp]
\begin{subfigure}{0.47\linewidth}
\centering
\includegraphics[width=\linewidth]{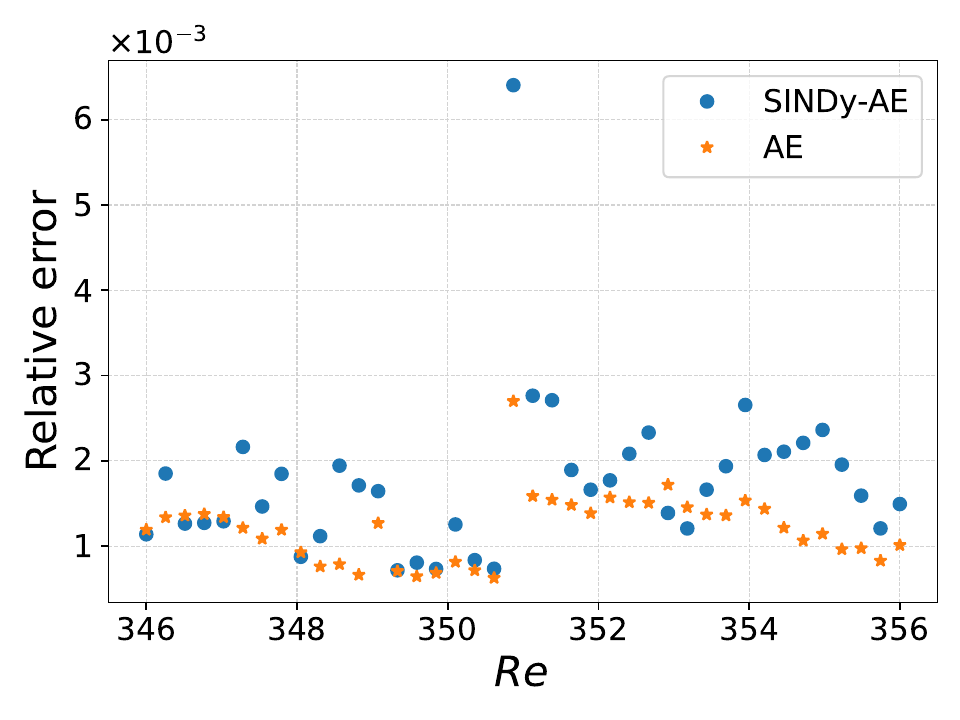}
\caption{$\bs{u}$}
\end{subfigure}
\hspace{0.6cm}
\begin{subfigure}{0.47\linewidth}
\centering
\includegraphics[width=\linewidth]{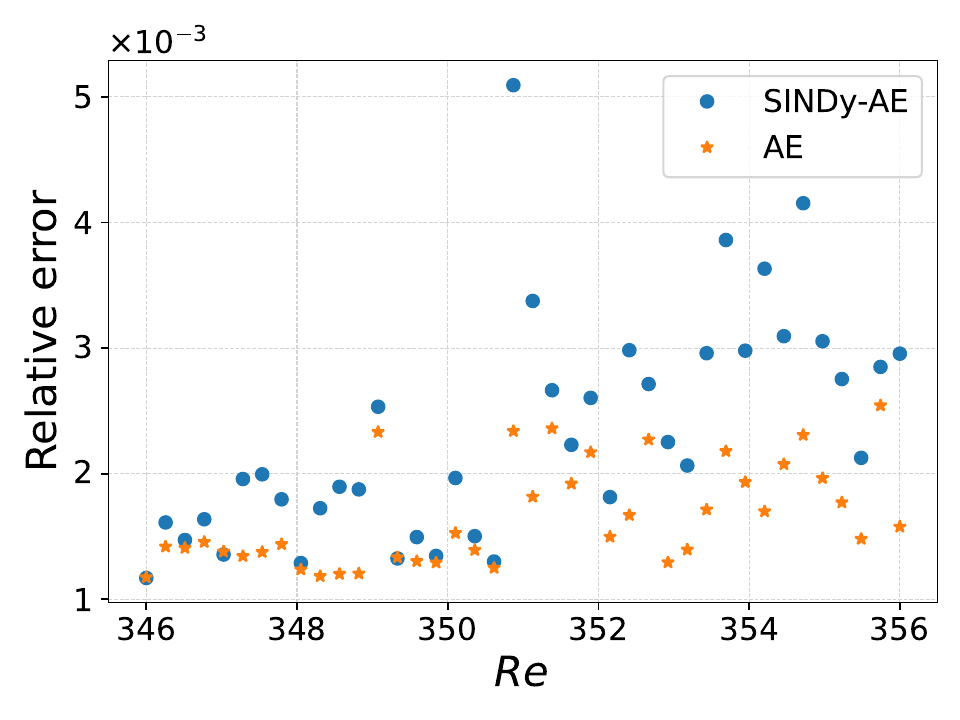}
\caption{$p$}
\end{subfigure}
\caption{Maximum error of the AE and SINDy-AE architecture with time integration on the POD coefficients of $\bs{u}$ and $p$ against $Re\in \mathbb{P}_h$.}
\label{figQUAINI_P_ERR}
\end{figure}

The comparison of the maximum errors introduced by the AE and SINDy-AE architectures are plotted against $Re\in\mathbb{P}_h$ in Figure \ref{figQUAINI_P_ERR}. We only observe a small discrepancy among the two, indicating that exploiting SINDy to identify and integrate the latent dynamics provides accurate results, with errors below $1\%$ for all values of $Re$, even when compared to the baseline AE reconstruction.

\begin{figure}[htbp]
\begin{subfigure}{0.47\linewidth}
\centering
\includegraphics[width=\linewidth]{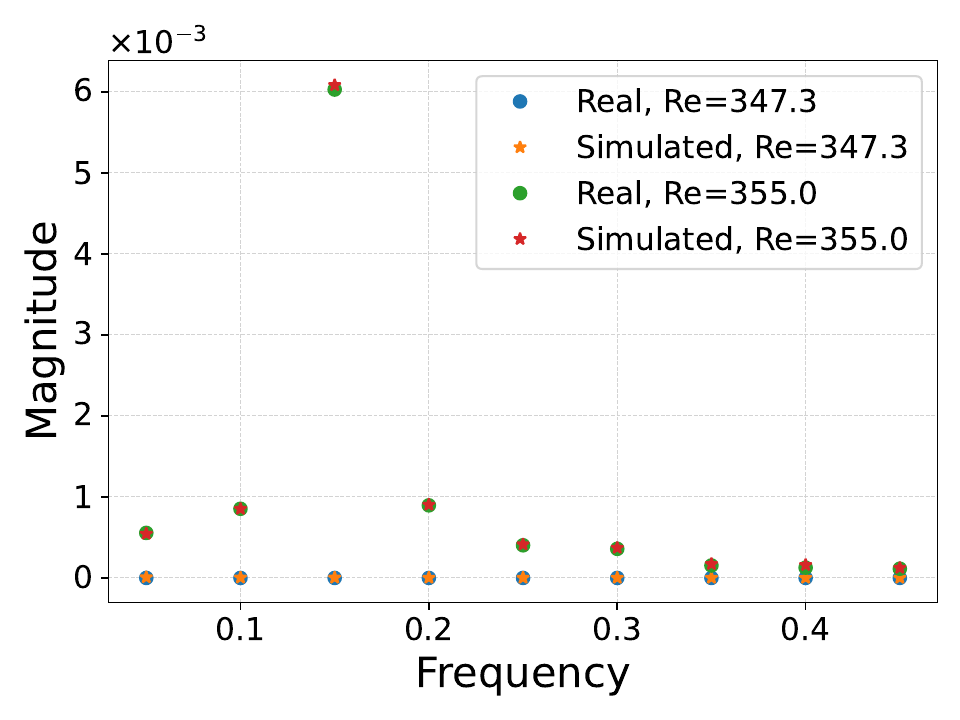}
\caption{$\bs{u}$}
\end{subfigure}
\hspace{0.6cm}
\begin{subfigure}{0.47\linewidth}
\centering
\includegraphics[width=\linewidth]{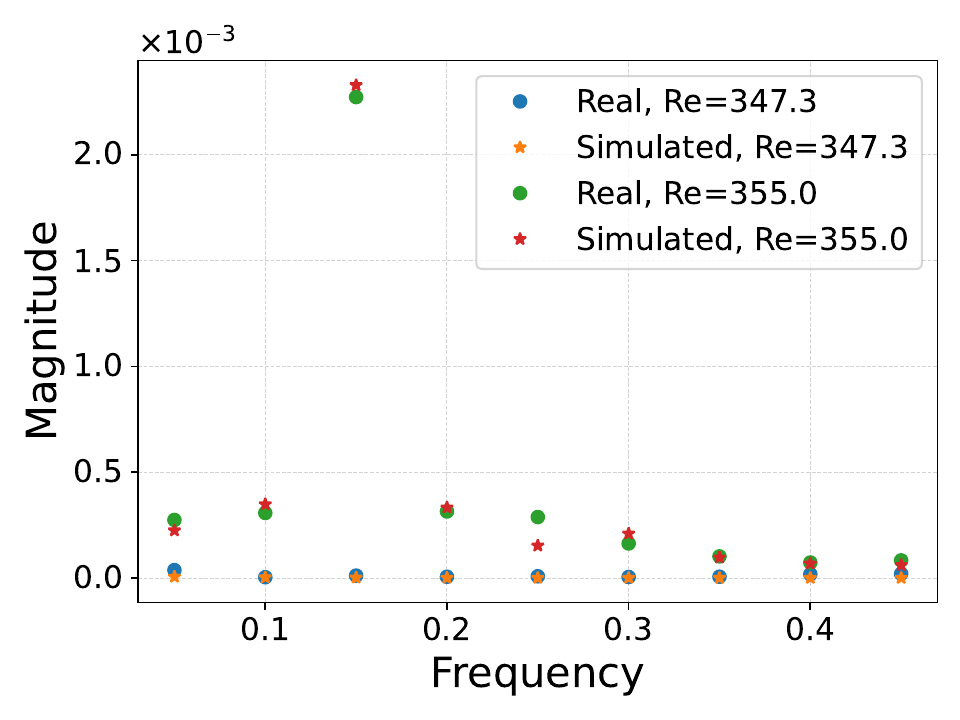}
\caption{$p$}
\end{subfigure}
\caption{Fourier analysis of the first POD coefficient of $\bs{u}$ and $p$ for $Re\in\mathbb{P}_h$.}
\label{fft_PARAM}
\end{figure}

In the same spirit of the frequency analysis depicted in Figure \ref{figFFT_energy}, we performed a Fourier analysis at the reduced level on the first POD coefficient for velocity and pressure fields. In Figure \ref{fft_PARAM} we show such analysis for the test $Re$ values, from which we can see both the accuracy between the real and simulated quantities, and the difference magnitude of the two parametric configurations, denoting the steady and unsteady behaviors.

\begin{figure}[htbp]
    \centering
    \includegraphics[width=0.5\textwidth]{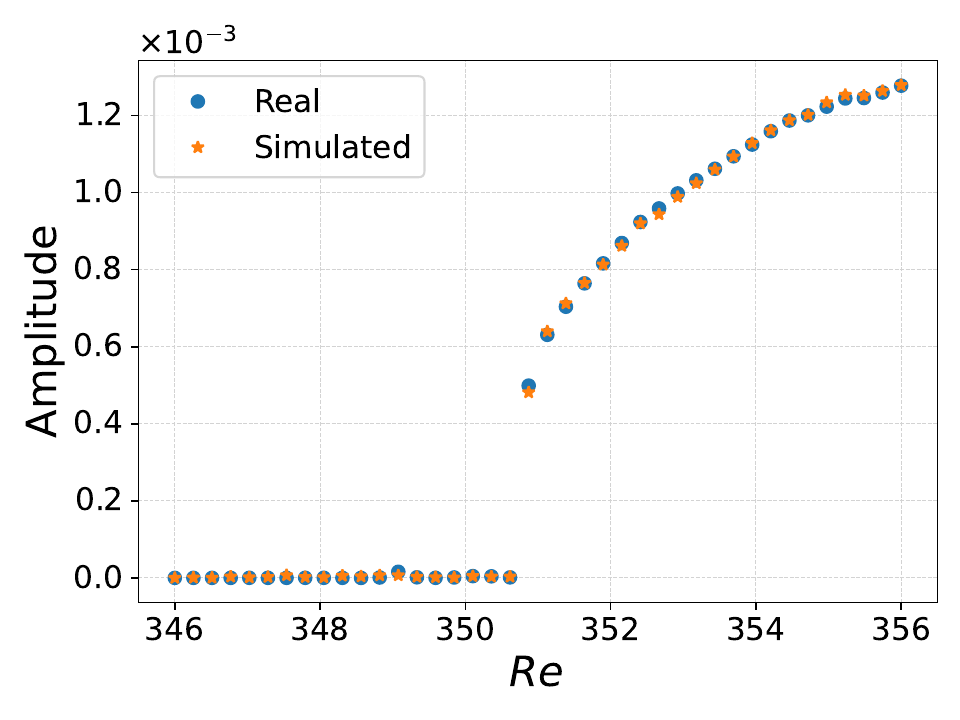}
    \caption{Real and simulated amplitudes against Reynolds number as defined in Equation \eqref{eqAMPLITUDE} for $Re\in\mathbb{P}_h$.}
    \label{figAMPLITUDE_PARAM}
\end{figure}

To conclude, Figure \ref{figAMPLITUDE_PARAM} shows the reconstruction of the bifurcation diagram using the amplitude defined in Equation \eqref{eqAMPLITUDE} as the quantity of interest, and identifying the bifurcation point as $Re^\ast\approx350.6$.

\subsubsection{SINDy-AE-POD architecture in the nonparametric case}
The identified ODEs for the previous systems, mixing both time and parameter dependency, have not been reported as they were not sparse enough and difficult to interpret.

This raises the question of whether this stems from an intrinsic difficulty in learning sparse representations for these bifurcating systems, or whether it is a limitation of the \gls{sindy} methodology when handling parameterized dynamics.
To proceed in this direction, we trained and tested a \gls{sindy}-AE model on a single trajectory.
We fixed a Reynolds number $Re=353.69$ in the bifurcating regime, and projected the snapshots of the corresponding dynamics onto the \gls{pod} coordinates.
For the \gls{sindy} algorithm, we used a library of first-order polynomials in $z_0$ and $z_1$, with the \gls{stlsq} threshold $\tau$ set to $1$.

To ensure complete system evolution, we analyzed the time interval $[1800, 2500]$, and used the first $90\%$ of the resulting dataset for training.
As before, since the time step between snapshots was $1$, we interpolated the POD coefficients using cubic natural spline functions evaluated at $dt=0.01$ second intervals.
In this case as well, we trained a separate \gls{sindy} model on latent variables and used ensembling for better robustness.

This way, we have been able to identify the following ODE:
\begin{equation*}
\begin{cases}
z_0' = -0.111 z_0 - 0.992 z_1 , \\
z_1' = 0.992 z_0 + 0.111 z_1 ,
\end{cases}
\end{equation*}
which exhibits a nice symmetry related to the periodic and oscillatory nature of the phenomena.

\begin{figure}[htbp]
    \centering
\includegraphics[width=0.9\textwidth]{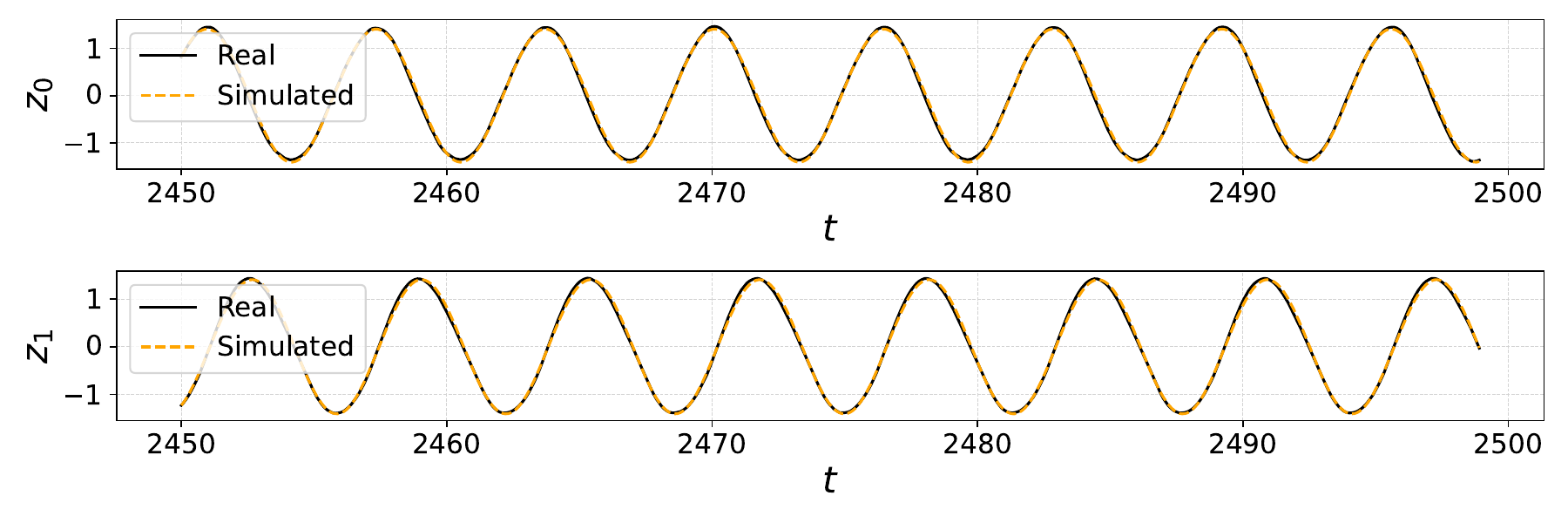}
    \caption{Latent trajectories on the integration window for $Re=353.69$.}
    \label{figQUAINI_LATENT_SINGLE}
\end{figure}

Moreover, we tested the model's ability to extrapolate in time beyond the training data by integrating over $[2450, 2500]$.
The great agreement between the real and simulated latent trajectories (below $1\%$), shown in Figure \ref{figQUAINI_LATENT_SINGLE}, confirms the possibility of capturing the periodic behavior in the latent dimension, performing excellently in terms of both recovered dynamics and error on the \gls{pod} coefficients.

Finally, being the identified ODE a sparse and analytically solvable system, it provides valuable insight into the underlying dynamics which could be exploited for more general frameworks.

\section{Conclusions and Perspectives}\label{sec:conclusions_and_perspectives}
In this work, we developed and investigated a \gls{rom} methodology for complex, time-dependent, nonlinear systems exhibiting bifurcation phenomena. Our main contribution is the successful application of the SINDy-AE-nested-POD architecture to challenging \gls{cfd} problems, demonstrating its effectiveness across different bifurcating scenarios and numerical discretization methods. In particular, we analyzed two distinct cases: a pitchfork bifurcation driven by the Coandă effect and a Hopf bifurcation for  channel flow problems, obtaining solutions using both FE and FV methods. The performance of our framework proved particularly strong, achieving robust identification of the pitchfork dynamics and excellent accuracy in capturing the Hopf behavior. Our investigation also provided novel insights into the challenges of parameterized system reduction, suggesting that the difficulties in obtaining effective latent representations through \gls{sindy} are inherent to the parameterized nature of the systems rather than limitations in trajectory representation.

During our study, several promising research directions emerged. A natural extension of this work would be to apply the SINDy-AE-nested-POD framework to even more challenging scenarios, such as compressible fluid dynamics or other fields characterized by complex dynamical behavior. Additionally, our experiments highlighted the high sensitivity of \gls{sindy} to hyperparameters. Recent advancements, including the nested \gls{sindy} method \cite{nested_SINDy} that integrates \gls{sindy} with neural networks, as well as techniques like WSINDy \cite{Messenger_2021} and LES-SINDy \cite{ZHENG2026114443}—which employ weak formulations or Laplace transforms to better manage discontinuities—offer enhanced robustness in the presence of noisy data. Although these methods incur higher computational costs, their potential synergy with suitable reduction techniques presents an exciting avenue for future research.

Another interesting research direction involves addressing the challenges associated with the investigation of parameterized systems, where some recent approaches have shown promise.
In particular, multi-objective optimization \cite{Lemus2024} appears to be effective in enabling meaningful identification of dynamics, even in the presence of noisy data and a limited number of trajectories.
This strategy could help mitigate a difficulty encountered in this work---specifically, enforcing sparsity at the latent level---thereby promoting interpretability of the identified dynamics in a parameterized setting.

Finally, a recent work by Conti, Kneifl et al.\ \cite{VVV, VVV_code} on a variational extension of \gls{sindy} (VINDy), combined with Variational Encoding of Noisy Inputs (VENI) and Variational Inference with Certainty Intervals (VICI), points to promising directions for dynamic system identification even in highly complex scenarios.

\section{Acknowledgments}
The authors thank Paolo Conti and Pasquale Claudio Africa for the fruitful discussions.
\textbf{FP} and \textbf{GR} acknowledge the support provided by the European Union-NextGenerationEU, in the framework of the iNEST-Interconnected Nord-Est Innovation Ecosystem (iNEST ECS00000043– CUP G93C22000610007) consortium and its CC5 Young Researchers initiative. The authors also like to acknowledge INdAM-GNCS and MIUR (Italian Ministry for University and Research) through FAREX-AROMA-CFD project, P.I.\ Prof. Gianluigi Rozza, for their support.

\renewcommand{\bibfont}{\small}
\printbibliography[]

\end{document}